\documentclass{amsart}
\usepackage{amssymb}
\usepackage{psfrag}
\usepackage{graphicx}
\usepackage{subfigure}
\usepackage{arydshln}
\usepackage{multirow}
\usepackage{xy}
\xyoption{all}

\newtheorem{defn}{Definition}[section]
\newtheorem{thm}[defn]{Theorem}

\newtheorem{lem}[defn]{Lemma}

\def\ker{\operatorname{Ker }}
\def\aut{\operatorname{Aut }}
\def\GL{\operatorname{GL }}

\def\B{\mathbf{B}}
\def\E{\mathbf{E}}
\def\Z{\mathbb{Z}}
\def\R{\mathbb{R}}
\def\P{\mathbf{P}}

\author{Byung Hee An}
\author{Ki Hyoung Ko}
\address{Department of Mathematics, Korea Advanced Institute
of Science and Technology, Daejeon, 305-701, Korea}
\email{\{anbyhee,knot\}@knot.kaist.ac.kr}
\thanks{This work was supported by the Korea Science and Engineering
Foundation(KOSEF) grant funded by the Korea government(MOST)
(No. R01-2006-000-10152-0)}

\begin{document}
\title{A family of representations of braid groups on surfaces}
\begin{abstract}
We propose a family of new representations of the braid groups on
surfaces that extend linear representations of the braid groups on
a disc such as the Burau representation and the Lawrence-Krammer-Bigelow
representation.
\end{abstract}

\maketitle

\section{Introduction and Preliminaries}

\subsection{Braid groups on surfaces}

In the 1920's, E.~Artin \cite{Art1} began to study braids as
mathematical objects and found a presentation of the group of braids
of $n$ strands under the binary operation concatenating
corresponding strands of two braids. Since then, braids have served
as convenient tools in various areas of science, especially in
physics. Fox and Neuwirth \cite{FoxN} generalized the braid groups
on arbitrary topological spaces $M$ via configuration spaces and
in this description, Artin's braid groups was on the plane $\mathbb R^2$.
Fadell and Neuwirth \cite{FN} found a fiber-bundle structure between
configuration spaces, and Birman \cite{Bir2} observed that there
is no interesting braid theory if the dimension of $M$ is
greater than $2$. It is possible to consider braid groups on
1-dimensional complex \cite{Gh}. But interesting braid groups are
often obtained on 2-dimensional manifolds and we are also interested
in these braid groups on surfaces in this arcticle.

Let $\Sigma(g,p)$ denote a compact, connected, orientable
2-dimensional manifolds of genus $g$ with $p$ boundary components.
We set $\Sigma= \Sigma(g,p)$. Let $\{z_1^0,\ldots, z_n^0\}$ be a set
of $n$ preferred distinct points in $\Sigma$ for $n\ge 0$ and  let
$\Sigma_n=\Sigma-\{z_1^0,\ldots, z_n^0\}$. We call $\Sigma_n$ the
surface $\Sigma$ {\em with $n$ punctures}.

For integers $n, k\ge 0$, we consider three types of {\em
configuration spaces} as follows:\\
The space of $k$-tuples of distinct points in $\Sigma_n$ denoted by
$$
P_{n,k}(\Sigma)=
\left\{(z_1,\cdots,z_k)\in \Sigma_n\times\cdots\times
\Sigma_n \mid z_i\neq z_j\ \mbox{for}\ i
\neq j\right\},
$$
the space of subsets of $k$ elements in $\Sigma_n$ denoted by
$$
B_{n,k}(\Sigma)=\left\{\{z_1,\cdots,z_k\}\subset \Sigma_n \right\},
$$
and the space $B_{n;k}(\Sigma)$ of pairs of disjoint subsets of $n$
elements and $k$ elements in $\Sigma$ denoted by
$$
B_{n;k}(\Sigma)=\left\{\left( \left\{ z_1,\dots,z_n \right\},\left\{z_{n+1},
\dots,z_{n+k} \right\} \right)\mid z_i\in \Sigma,\: z_i\ne z_j\mbox{ for }
i\ne j\right\}.
$$
It is easy to see that
$B_{n,k}(\Sigma)=P_{n,k}(\Sigma)/\mathbf{S}_k$ and
$B_{n;k}(\Sigma)=P_{0,n+k}(\Sigma)/\mathbf{S}_n\times\mathbf{S}_k$
where the symmetric group $\mathbf{S}_k$ acts on $P_{n,k}(\Sigma)$
by permuting components of a $k$-tuple and similarly
$\mathbf{S}_n\times\mathbf{S}_k\subset \mathbf S_{n+k}$ acts on
$P_{0,n+k}(\Sigma)$.

The braid groups on a surface $\Sigma$ are defined by the
fundamental groups of configuration spaces. Choose a basepoint
$\{z_{n+1}^0,\ldots,z_{n+k}^0\}$ in $\partial \Sigma$ if
$\partial\Sigma\neq\emptyset$. If $\Sigma$ is closed, then place it
anywhere in $\Sigma_n$. The {\em pure $k$-braid group on $\Sigma_n$}
is defined and denoted by
$$\P_{n,k}(\Sigma)=\pi_1\left(P_{n,k}(\Sigma),(z_{n+1}^0,
\ldots,z_{n+k}^0)\right).$$
Similarly, the {\em (full) $k$-braid group on $\Sigma_n$} is given
by
$$\B_{n,k}(\Sigma)=\pi_1\left(B_{n,k}(\Sigma),\{z_{n+1}^0,
\ldots,z_{n+k}^0\}\right),$$
and the {\em intertwining $(n,k)$-braid group on $\Sigma$} is given
by
$$\B_{n;k}(\Sigma)=\pi_1\left(B_{n;k}(\Sigma),\left(\{z_1^0,\ldots,z_n^0\},
\{z_{n+1}^0,\ldots, z_{n+k}^0\}\right)\right).$$ It is sometimes
easier to understand if these groups are regarded as subgroups of
$\B_{0,n+k}(\Sigma)$. The intertwining $(n,k)$-braid group
$\B_{n;k}(\Sigma)$ is the preimage of $\mathbf S_n\times \mathbf
S_k$ under the canonical projection$:\B_{0,n+k}(\Sigma)\to \mathbf
S_{n+k}$. In addition, $\B_{n,k}(\Sigma)$ is the subgroup of
$(n+k)$-braids in $\B_{n;k}(\Sigma)$ that become trivial by
forgetting the last $k$ strands and $\P_{n,k}(\Sigma)$ is the
subgroup of $(n+k)$-braids in $\B_{n,k}(\Sigma)$ that are pure, that
is, the induced permutation is trivial. If the surface $\Sigma$ is
the 2-disc $D$ 
we will call the braid groups {\em classical} . For example, $\B_{0,n}(D)$
denotes the classical $n$-braid group studied by E.~Artin.

In the 60's and 70's, presentations for braid groups on various
surfaces are found: on the 2-sphere and the projective plane by
Fadell and Van Buskirk \cite{FV, Van}, on the torus by Birman
\cite{Bir2}, on all closed surfaces by Scott \cite{Sco}. The study
on braid groups on surfaces has been revived recently.
Gonz\'alez-Meneses \cite{Gon} found new presentations of the braid
groups on surfaces and Bellingeri \cite{Bel, BG} found positive
presentations of the braid groups $\B_{n,k}(\Sigma)$ for all
surfaces $\Sigma$ with or without boundaries. In this article, we
are interested in braid groups that contain the classical braid
groups as subgroups. So we will consider surfaces with nonempty
boundary and will use Bellingeri's presentations.

We remark that boundary components of a surface can be traded with
punctures when we consider braid groups. Let $\Sigma=\Sigma(g,p)$
and $\Sigma'=\Sigma(g,p+q)$. Then there are continuous maps
$i:\Sigma_q\to\Sigma'$ and $j:\Sigma'\to\Sigma_q$ that are homotopy
inverses each other. The induced maps $\bar i:B_{n+q,k}(\Sigma)\to
B_{n,k}(\Sigma')$ and $\bar j:B_{n,k}(\Sigma')\to B_{n+q,k}(\Sigma)$
on configuration spaces are also homotopy inverses each other and
induce isomorphisms $\bar i_*, \bar j_*$ on fundamental groups
\cite{Bel, PR}. Therefore we may assume $\Sigma=\Sigma(g,1)$ by
treating all but one boundary component as punctures whenever we
deal with a surface with nonempty boundary.

Throughout the article, we use Bellingeri's presentation \cite{Bel}
for the braid group $\B_{n,k}(\Sigma(g,1))$ given as follows:
\begin{itemize}
\item Generators :
$\sigma_1,\dots,\sigma_{k-1},a_1,\dots,a_g,b_1,\dots,b_g,\zeta_1,\dots,
\zeta_n$.
\item Relations
\begin{itemize}
\item [BR1:] $\left[\sigma_i,\sigma_j\right]$,\quad $|i-j|\ge 2$
\item [BR2:] $\sigma_i\sigma_j\sigma_i=\sigma_j\sigma_i\sigma_j$,
\quad $|i-j|=1$
\item [CR1:] $\left[a_r,\sigma_i\right]$,
$\left[b_r,\sigma_i\right]$,
$\left[\zeta_t,\sigma_i\right]$, \quad $i>1$
\item [CR2:] $\left[a_r,\sigma_1 a_r\sigma_1\right]$,
$\left[b_r,\sigma_1 b_r\sigma_1\right]$,
$\left[\zeta_t,\sigma_1 \zeta_t\sigma_1\right]$
\item [CR3:] $[a_r,\sigma_1^{-1} a_s\sigma_1]$,
$[a_r,\sigma_1^{-1} b_s\sigma_1]$, $[b_r,\sigma_1^{-1}
a_s\sigma_1]$, $[b_r,\sigma_1^{-1} b_s\sigma_1]$, \quad
$r<s$
\item [] $[a_r,\sigma_1^{-1} \zeta_u\sigma_1]$,
$[b_r,\sigma_1^{-1} \zeta_u\sigma_1]$,
$[\zeta_t,\sigma_1^{-1} \zeta_u\sigma_1]$, \quad $t<u$
\item [SCR:] $\sigma_1 b_r\sigma_1 a_r\sigma_1=a_r\sigma_1 b_r$.
\end{itemize}
\end{itemize}

The corresponding result for configuration spaces of pure braids by
Fadell and Neuwirth can be generalized to show that the projection
$B_{n;k}(\Sigma)\to B_{0,n}(\Sigma)$ is a fiber bundle with fiber
$B_{n,k}(\Sigma)$. Except for $S^2$ and $\R P^2$, Guaschi and Gon\c
calves \cite{GG1} completely determined when the short exact
sequences of braid groups derived from Fadell-Neuwirth fibrations
split. In particular, the short exact sequence derived from the
fibration above
$$1\to\B_{n,k}(\Sigma)\to\B_{n;k}(\Sigma)\to\B_{0,n}(\Sigma)\to1$$
always splits for $k\ge 1$ if $\Sigma$ has nonempty boundaries.

The braid groups are closely related to the mapping class groups.
Birman \cite{Bir} determined when surface braid group embeds into
the corresponding mapping class group. In particular, if
$\partial\Sigma\neq\emptyset$, $\B_{0,n}(\Sigma)$ embeds into the
mapping class group on $\Sigma_n$ and so an $n$-braid on $\Sigma$
can be regarded as a homeomorphism of $\Sigma$ that preserves the
set of $n$ punctures.

\subsection{Representations of braid groups}
The classical braid groups have various representations that can be
as simple as taking exponent sums or taking induced permutations.
The braid action on the punctured disk $D_n$ gives rise a faithful
representation into automorphism groups of free groups and a
characterization of automorphisms coming from braid actions is
possible. Each representation serves its own purpose. It is common
to try to construct a linear representation to have a better
understanding of a given group via matrices over a certain commutative ring and
their multiplications.

For the classical braid groups, linear representations are abundant.
Burau in 1936 and Gassner in 1961 discovered linear representations
of $\B_{0,n}(D)$ and $\P_{0,n}(D)$, respectively. These representations
are derived from braid actions on homologies of appropriate
coverings of $D_n$. These representations take $(n-1)\times(n-1)$
matrices that can also be computed via Fox's free differential
calculus on automorphisms of free groups mentioned above. The Burau
representation is faithful for $n\le 3$ but not faithful for $n\ge
5$ \cite{Big1}. The faithfulness of Gassner representation is known
only for $n\le 3$.

Lawrence \cite{Law} discovered a family of linear representations of
$\B_{0,n}(D)$ via a monodromy on a vector bundle over $P_{n,k}(D)$.
Krammer \cite{Kra1} defined a free $\Z[q^{\pm1},t^{\pm1}]$-module
$V$ using {\em forks} and relations between them and he proved that
the braid group acts on $V$ faithfully for braid index $4$ using
algebraic and combinatorial argument. This representation is
essentially the same as the one considered by Lawrence for $k=2$ but
uses configuration space $B_{n,2}(D)$ instead of $P_{n,2}(D)$ and is now
called the {\em Lawrence-Krammer-Bigelow representation}. Bigelow
reinterpreted this representation using covering spaces and covering
transformation groups instead of vector bundles and local
coefficients. Then the monodromy corresponds to the braid action on
homology groups of covering spaces as the Burau representation and
the Gassner representation were obtained. Bigelow \cite{Big2}
constructed a linear representation using homology group $H_2(\tilde
B_{n,k}(D))$ of the covering space $\tilde B_{n,k}(D)$ whose covering
transformation group is $\langle q\rangle\oplus\langle t\rangle$,
and he proved that $\R\otimes V$ is isomorphic to $\R\otimes
H_2(\tilde B_{n,k}(D))$. Moreover, Krammer \cite{Kra2} and Bigelow
\cite{Big2} independently proved that Lawrence-Krammer-Bigelow
representation is faithful for all $n\ge 1$, and so the classical
braid groups are linear. Furthermore, Bigelow and Budney \cite{BB}
proved that the mapping class group of genus 2 surface has a
faithful linear representation using the Lawrence-Krammer-Bigelow
representation and suitable branched covering.
However, Paoluzzi and Paris showed that there is a difference
between $V$ and $H_2(\tilde B_{n,k})(D)$ as
$\Z[q^{\pm1},t^{\pm1}]$-module for $n\ge 3$ and found basis for
$\Z[q^{\pm1},t^{\pm1}]$-module $H_2(\tilde B_{n,k}(D))$ and so the
exact definition of ``Lawrence-Krammer-Bigelow representation" became little
ambiguous.

For any $k\ge 1$, Bigelow in \cite{Big4} considered the braid
action on the Borel-Moore homology
group $H_k^{BM}(\tilde B_{n,k}(D))$. He obtained a family
of representations via the induced action on the image of
$H_k^{BM}(\tilde B_{n,k}(D))$ in  $H_k^{BM}(\tilde B_{n,k}(D), \partial\tilde B_{n,k}(D))$.
For simplicity, we will directly
consider the braid action on the free module $H_k^{BM}(\tilde
B_{n,k}(D))$ whose basis can be easily described by forks to obtain a
linear representation and we will call these representations {\em
homology linear representations}. The Burau representation and the
Lawrence-Krammer-Bigelow representation of $\B_{0,n}(D)$ are homology
linear representations $$\Phi_k:\B_{0,n}(D)\to
\GL\left(\left(\begin{matrix} n+k-2\\k
\end{matrix}\right),\Z[q^{\pm1},t^{\pm1}]\right)$$ obtained from
braid action on homologies of covers of $B_{n,k}(D)$ when $k=1$ ($t=1$
in this case) and $k=2$, respectively \cite{Big4}.

In this article, we construct a family of representations of braid
groups on surface with nonempty boundary that extend homological
linear representations of the classical braid group. In \S2, we first
try to follow the way how the homology linear representations of
$\B_{0,n}(D)$ was constructed via a covering of the configuration space
$B_{n,k}(D)$. In the case of the disk, the braid action automatically
commutes with covering transformation or braids act trivially on
local coefficients in other words. However in case of surfaces of
genus $\ge 1$, this condition forces the variable $q$ equal $1$ (see Lemma~\ref{lem:doublestaristoostrong}).
Then the braid action becomes almost trivial. For example, if $k=1$,
the action of $\sigma_i^2$ is trivial. In order to get around this
problem, we introduce the intertwining braid group
$\B_{n;k}(\Sigma)$ in \S3 to replace  $\B_{n,k}(\Sigma)$. As we
mentioned earlier, this group is a semidirect product of
$\B_{n,k}(\Sigma)$ and $\B_{0,n}(\Sigma)$.  Although the braid action does not preserve the
local coefficient given by the $\B_{n,k}(\Sigma)$ factor, the
$\B_{0,n}(\Sigma)$ factor of $\B_{n;k}(\Sigma)$ can adjust the
coefficient so that the braid action becomes compatible.
The coefficient ring for
homology representations will also be extended in order to give a room to
control the braid action in the expense of giving up its commutativity
so that it becomes more interesting and
still preserves coefficients. Eventually we
obtain representations of braid groups on surfaces that extend
homology linear representations of the classical braid group (see Theorem~\ref{thm:representation}).
And we explicitly compute the representations in the form of
matrices using a geometric argument. We extend the intersection
pairing between $H_k^{BM}(\tilde B_{n,k}(D))$ and its dual space
$H_k(\tilde B_{n,k}(D), \partial \tilde B_{n,k}(D))$ and use bases for the
two spaces that are described by ``forks" and ``noodles" (see Theorem~\ref{thm:computation}).

In \S4, we argue that the construction of our
representations is natural and there are no other alternatives if
one wants to obtain an extension of the homological representation
using covers of the configuration space $B_{n,k}(D)$.
We show that the
intertwining braid group $\B_{n;k}(\Sigma)$ is the normalizer of
$\B_{n,k}(\Sigma)$ in $\B_{0,n+k}(\Sigma)$ so that the
intertwining braid group $\B_{n;k}(\Sigma)$ is such a group that is
unique and maximal if we ignore a meaningless extension (see Theorem~\ref{thm:uniqueness}).
The coefficient ring for our representations is
the integral group ring of a quotient group of $\B_{n;k}(\Sigma)$.
For $k\ge 3$, the quotient group is uniquely determined in order to extend homology linear representations of the classical braid group.
For $k=1, 2$, we use the simplest one as long as it serves
our purpose. Then we show that the braid action on the quotient group
is virtually unique (see Theorem~\ref{thm:character}).

Our construction involving the group extension $\B_{n;k}(\Sigma)$ of
$\B_{n,k}(\Sigma)$ is purely algebraic and does not carry a good
geometric interpretation. As a result, some of useful geometric
tools are not available. For example. the intersection pairing
mentioned above is not invariant under the braid group action. This
seems to make it difficult to discuss properties of our
representations such as faithfulness and irreducibility. Although
the corresponding representation of the classical braid group is
faithful for $k=2$ and irreducible for $k\le2$ \cite{Jon,Zin}, both
faithfulness and irreducibility for our representations are beyond
the scope of the current article.

\section{Homology linear representations}

We first review the construction of homology linear representations
of the classical braid group $\B_{0,n}(D)$ using the configuration
space $B_{n,k}(D)$ and then discuss the difficulty in extending these
homology linear representations to the braid group
$\B_{0,n}(\Sigma)$ on a surface $\Sigma$ with nonempty boundary. As
we noted earlier, boundary components can be traded with punctures.
From now on, we assume that $\Sigma$ denotes a compact, connected,
oriented surface with exactly one boundary component and that $n$
and $k$ are positive integers.

\subsection{Homology linear representations of classical braid group}
Let $\phi:\B_{n,k}(D)\to G$ be an epimorphism onto a
group $G$. Consider the covering $p:\tilde B_{n,k}(D)\to
B_{n,k}(D)$ corresponding to $\ker{\phi}$. Since the classical braid
group embeds into the mapping class group of the punctured disk
$D_n$, we may assume we have a homeomorphism $\bar\beta:B_{n,k}(D)\to
B_{n,k}(D)$ for each $\beta\in\B_{0,n}$. By the lifting criteria,
$\bar\beta$ lifts to $\tilde\beta:\tilde B_{n,k}(D)\to\tilde B_{n,k}(D)$ if
and only if $\bar\beta_*(\ker\phi)\subset \ker\phi$. This is equivalent that
there is an induced automorphism $\beta_\sharp$ on $G$ such that
$\beta_\sharp\phi=\phi\bar\beta_*$.

Now we consider {\em Borel-Moore homology} \cite{BM,HR}
defined by
$$H_\ell^{BM}(\tilde B_{n,k}(D))=\lim_{\leftarrow}H_\ell(\tilde B_{n,k}(D),
p^{-1}(B_{n,k}(D)\setminus A))$$
where the inverse limit is taken over all compact subsets $A$ of
$B_{n,k}(D)$.

The middle-dimensional homology group $H_k^{BM}(\tilde B_{n,k}(D))$ is
a free $\Z[G]$-module of rank $\binom{n+k-2}{k}$ (see \cite{Big4})
and $\tilde\beta$ induces a map $\tilde\beta_*:H_k^{BM}(\tilde
B_{n,k}(D))\to H_k^{BM}(\tilde B_{n,k}(D))$ such that
$$\tilde\beta_*(y c)=\beta_\sharp(y)\tilde\beta_*(c)$$
for $y\in G$ and $c\in H_k^{BM}(\tilde B_{n,k}(D))$. Thus the map
$\tilde\beta_*$ is a $\Z[G]$-module homomorphism if and only if
$\beta_\sharp(y)=y$ for all $y\in G$
if and only if
\begin{equation}
\phi=\phi\bar\beta_*\quad \mbox{for all }\beta\in \B_{0,n}\tag{$*$}.
\end{equation}

Notice that the condition ($*$) also implies
$\bar\beta_*(\ker\phi)\subset \ker\phi$. Here we need to know that the
induced homomorphism $\tilde\beta_*$ depends only on the isotopy
class of the homeomorphism $\beta$. In fact, since $D$ has a
boundary, we choose the basepoint $\left\{
z_{n+1}^0,\cdots,z_{n+k}^0 \right\}$ of $B_{n,k}(D)$ in $\partial D$
and then the isotopy preserves the basepoint and gives the same
induced map $\tilde\beta_*$. Consequently, if we choose a group $G$
and an epimorphism $\phi:\B_{n,k}(D)\to G$ satisfying ($*$), we obtain
a family of representation
$\Phi_k:\B_{0,n}(D)\to\aut_{\Z[G]}\left(H_k^{BM}(\tilde
B_{n,k}(D))\right)$ into the group of $\Z[G]$-module automorphisms on
$H_k^{BM}(\tilde B_{n,k}(D))$ defined by
$$\Phi_k(\beta)=\tilde{\beta}_*:H_k^{BM}(\tilde B_{n,k}(D))\to
H_k^{BM}(\tilde B_{n,k}(D)).$$

Because we want to obtain a linear representation, the group $G$
should be abelian. By the presentation given in \S1.1, $\B_{n,k}(D)$ is
generated by $\zeta_1,\ldots,\zeta_n$,
$\sigma_1,\ldots,\sigma_{k-1}$.
Suppose that $\phi:\B_{n,k}(D)\to G$ be an epimorphism satisfying the condition $(*)$ onto an abelian group $G$.
Each generator $\sigma_i$ of $\B_{0,n}(D)$ acts trivially on
$\B_{n,k}(D)$ except
$$
    (\bar\sigma_i)_*(\zeta_i)=\zeta_i \zeta_{i+1}\zeta_i^{-1},
    (\bar\sigma_i)_*(\zeta_{i+1})=\zeta_i.
$$

Then the condition $(*)$ implies that $\phi(\zeta_i)=\phi((\bar\sigma_i)_*(\zeta_{i+1}))=\phi(\zeta_{i+1})$.
Hence for $k=1$,  $G$ is a quotient of $\langle q\rangle$ and
$\phi(\zeta_i)=q$ for $i=1,\ldots,n$.
For $k\ge 2$, $G$ is a quotient of $\langle
q\rangle\oplus\langle t\rangle$ and $\phi(\zeta_i)=q,
\phi(\sigma_j)=t$ for $i=1,\ldots,n$ and $j=1,\ldots,k-1$.

We define a group $G_D$ and an epimorphism $\phi_D:\B_{n,k}(D)\to G_D$ depending only on $k$ as follows:
$$
\phi_D:\B_{n,k}(D)\to G_D=\begin{cases}
\langle q\rangle& k=1\\
\langle q\rangle\oplus\langle t\rangle& k\ge 2
\end{cases}.
$$

\begin{thm}\cite{Big4, Law}
    \label{thm:representationondisc}
Let $\phi_D:\B_{n,k}(D)\to G_D$ be the epimorphism defined as above.
Then there is a homomorphism
$$
\Phi_k:\B_{0,n}(D)\to \aut_{\Z[G_D]}\left(H_k^{BM}(\tilde B_{n,k}(D))\right).
$$
In fact, $\Phi_1$ is the Burau representation and $\Phi_2$ is the
Lawrence-Krammer-Bigelow representation.
\end{thm}

\subsection{Naive extension to braid groups on surfaces}

Let $\Sigma$ be a surface of genus $g\ge 1$ with nonempty boundary.
Besides the two reasons mentioned at the end of \S1.1, there are one
more reason that makes the assumption $\partial\Sigma\ne\emptyset$
necessary. Suppose that $\partial\Sigma=\emptyset$ and $\beta\in
\B_{0,n}(\Sigma)$ uniquely determines the isotopy class of a
homeomorphism $\bar\beta:B_{n,k}(\Sigma)\to B_{n,k}(\Sigma)$. Then we
must choose the basepoint $\{z_{n+1}^0,\dots,z_{n+k}^0\}$ in the
interior of $\Sigma$. We can easily find a homeomorphism
$\bar\beta:B_{n,k}(\Sigma)\to B_{n,k}(\Sigma)$ that is isotopic to the
identity via an isotopy that does not preserve the basepoint. Then
$\beta$ represents the identity element in $\B_{0,n}(\Sigma)$ but
$\tilde\beta_*:H_{k}^{BM}(\tilde B_{n,k}(\Sigma))\to
H_{k}^{BM}(\tilde B_{n,k}(\Sigma))$ may be nontrivial. Thus no
representation can be obtained in this way if
$\partial\Sigma=\emptyset$.

We need to define when we say that a representation of the braid
group $\B_{0,n}(\Sigma)$ extends homology linear representations of
the classical braid groups.

\begin{defn}
\label{defn:extension}
Given a ring $R$, let $M$ be an
$R$-module on which the braid group $\B_{0,n}(\Sigma)$ acts as
$R$-module isomorphisms.
The $R$-module $M$ is an {\em extension} of
homology linear representations of the classical braid groups
$\B_{0,n}(D)$ if there exists a $\Z[G_D]$-submodule $M'$ of $M$ such that
\begin{itemize}
\item[(i)] It is invariant under the action by the subgroup $\B_{0,n}(D)$ of
$\B_{0,n}(\Sigma)$;
\item[(ii)]  For some $k\ge 1$,
$R$ contains $\Z[G_D]$ as a subring and
there is a $\Z[G_D]$-isomorphism from
$H_k^{BM}(\tilde B_{n,k}(D))$ to $M'$ that commutes with the
$\B_{0,n}(D)$ action.
\end{itemize}
\end{defn}

As in the classical braid cases, we have to look at the action of $\B_{0,n}(\Sigma)$ on $\B_{n,k}(\Sigma)$. The following lemma helps us to observe the action we desire.
\begin{lem}\cite{Bir, FN, GG1}
    \label{lem:splits}
Let $\pi_n:B_{n;k}(\Sigma)\to B_{0,n}(\Sigma)$ be the projection
onto the first $n$ coordinates. Then the space $B_{n;k}(\Sigma)$ is
a fiber bundle with fiber $B_{n,k}(\Sigma)$ and the induced short
exact sequence
$$
\xymatrix{
    1\ar[r]&\B_{n,k}(\Sigma)\ar[r]&\B_{n;k}(\Sigma)
    \ar[r]^{(\pi_n)_*}&\B_{0,n}(\Sigma)\ar[r]&1
    }
$$
splits for all $k\ge 1$.
\end{lem}

This lemma shows that how to decompose a braid
$\beta\in\B_{n;k}(\Sigma)$ into a product $\beta=\beta_1\beta_2$ for
$\beta_1\in\B_{0,n}(\Sigma)$ and $\beta_2\in\B_{n,k}(\Sigma)$.
Let $\iota:\B_{0,n}(\Sigma)\to\B_{n;k}(\Sigma)$ be the splitting map.
Then the above lemma shows that $\B_{n;k}(\Sigma)$ can be generated by the following two sets:
$$
X_1=\{\bar\sigma_1,\dots,\bar\sigma_{n-1},\bar
\mu_1,\dots,\bar \mu_{g}, \bar \lambda_1,\dots,\bar \lambda_g\},
$$
$$
X_2=\{\sigma_1,\dots,\sigma_{k-1},\zeta_1,\dots,
\zeta_n,\mu_1,\dots,\mu_{g},\lambda_1,\dots,\lambda_g\}
$$
where the generators in $X_1$ are the images of generators in $\B_{0,n}(\Sigma)$ under the inclusion map $\iota$.

Then the action of $\B_{0,n}(\Sigma)$ on $\B_{n,k}(\Sigma)$ is equivalent to the conjugate action in $\B_{n;k}(\Sigma)$ if we regard these two groups as subgroups of $\B_{n;k}(\Sigma)$.
The following lemma shows how
$\B_{0,n}(\Sigma)$ acts on $\B_{n,k}(\Sigma)$ and its proof is
straightforward and omitted.
\begin{lem}
\label{lem:braidaction} Each generator of $\B_{0,n}(\Sigma)$ acts on
$\B_{n,k}(\Sigma)$ as follows.
\begin{itemize}
\item [(1)] For $1\le i\le n-1$,

$
\begin{array}{rll}
(\bar\sigma_i)_*(\zeta_t)&=&
\left\{
\begin{array}{ll}
    \zeta_i\zeta_{i+1}\zeta_i^{-1}&t=i\\
\zeta_{i}&t=i+1
\end{array}
\right.
\end{array}
$
\item [(2)] For $1\le r\le g$,

$
\begin{array}{rll}
(\bar \mu_r)_*(\zeta_1)&=&\mu_r\zeta_1\mu_r^{-1}\\
(\bar \mu_r)_*(\mu_s)&=& \left\{
\begin{array}{ll}
\mu_r\zeta_1 \mu_r\zeta_1^{-1}\mu_r^{-1}&s=r\\
\left[\mu_r,\zeta_1\right]\mu_s\left[\mu_r,\zeta_1\right]^{-1}& r<s
\end{array}
\right.\\
(\bar \mu_r)_*(\lambda_s)&=& \left\{
\begin{array}{ll}
\lambda_r\mu_r\zeta_1^{-1}\mu_r^{-1}&s=r\\
\left[\mu_r,\zeta_1\right]\lambda_s\left[ \mu_r,\zeta_1 \right]^{-1}& r<s
\end{array}
\right.
\end{array}
$
\item [(3)] For $1\le r\le g$,

$
\begin{array}{rll}
(\bar \lambda_r)_*(\zeta_1)&=&\lambda_r\zeta_1\lambda_r^{-1}\\
(\bar \lambda_r)_*(\mu_s)&=& \left\{
\begin{array}{ll}
\lambda_r \zeta_1 \lambda_r^{-1} \mu_r \zeta_1 \lambda_r \zeta_1^{-1} \lambda_r^{-1}&s=r\\
\left[\lambda_r,\zeta_1\right]\mu_s\left[ \lambda_r,\zeta_1 \right]^{-1}& r<s
\end{array}
\right.\\
(\bar \lambda_r)_*(\lambda_s)&=& \left\{
\begin{array}{ll}
\lambda_r\zeta_1\lambda_r\zeta_1^{-1}\lambda_r^{-1}&s=r\\
\left[\lambda_r,\zeta_1\right]\lambda_s\left[ \lambda_r,\zeta_1 \right]^{-1}& r<s
\end{array}
\right.
\end{array}
$
\item [(4)] All other generators act trivially.
\end{itemize}
\end{lem}

We can find the presentation for $\B_{n;k}(\Sigma)$ using the above lemma as follows.
\begin{lem}
The braid group $\B_{n;k}(\Sigma)$ admits the following presentation
\begin{itemize}
\item Generators:
\begin{itemize}
\item[]
$X_1=\{\bar\sigma_1,\dots,\bar\sigma_{n-1},\bar
\mu_1,\dots,\bar \mu_g, \bar \lambda_1,\dots,\bar \lambda_g\}$
\item[]
$X_2=\{\sigma_1,\dots,\sigma_{k-1},\zeta_1,\dots,
\zeta_n,\mu_1,\dots,\mu_g,\lambda_1,\dots,\lambda_g\}$
\end{itemize}
\item Relations:
\begin{itemize}
\item [(i)] BR1 through SCR among generators in $X_1$
\item [(ii)] BR1 through SCR among generators in $X_2$
\item [(iii)] $ \bar x^{-1}y \bar x= (\bar x_*)(y)$ for all $\bar x\in X_1, y\in X_2$
\end{itemize}
\end{itemize}
where the relations BR1 through SCR are from Bellingeri's
presentation in \S1.1 and the action by $\bar x_*$ is given in Lemma~\ref{lem:braidaction}.
\label{lem:presentationofb_n;k}
\end{lem}
\begin{proof}
By Lemma~\ref{lem:splits}, the intertwining braid group
$\B_{n;k}(\Sigma)$ is a semidirect product of the normal subgroup
$\B_{n,k}(\Sigma)$ and $\B_{0,n}(\Sigma)$ where $\B_{0,n}(\Sigma)$
acts on $\B_{n,k}(\Sigma)$ by conjugation as shown in Lemma~\ref{lem:braidaction}.
Then it is easy to show that the semidirect
product $\B_{n;k}(\Sigma)$ admits the desired presentation. 
\end{proof}

For surfaces, the condition $(*)$ implies an undesirable consequence
as shown in the following lemma.
\begin{lem}
\label{lem:doublestaristoostrong} Let $\phi:\B_{n,k}(\Sigma)\to G$
be an epimorphism satisfying $\phi=\phi\bar\beta_*$  for any
$\beta\in\B_{n,k}(\Sigma)$. Then for $i=1,\ldots,n$,
$$\phi(\zeta_i)=1.$$
\end{lem}
\begin{proof}
As seen earlier, the hypothesis on $\phi$ implies that  for all
$y\in G$ and $r=1,\ldots,g$,
$$(\mu_r)_\sharp(y)=y.$$
But by Lemma~\ref{lem:braidaction} (2), we have
\begin{align*}
(\mu_r)_\sharp\phi(\lambda_r)&=\phi((\bar \mu_r)_*(\lambda_r))
=\phi(\lambda_r\mu_r\zeta_1^{-1}\mu_r^{-1})\\
&=\phi(\lambda_r)\phi((\bar \mu_r)_*(\zeta_1^{-1}))
=\phi(\lambda_r)(\mu_r)_\sharp(\phi(\zeta_1^{-1}))
=\phi(\lambda_r)\phi(\zeta_1^{-1}).
\end{align*}
Since $(\mu_r)_\sharp(\phi(\lambda_r))=\phi(\lambda_r)$ by hypothesis,
$\phi(\zeta_1)=1$ and so $\phi(\zeta_i)=1$ for all $1\le i\le n$ by
Lemma~\ref{lem:braidaction} (1).
\end{proof}

The previous lemma says that the condition $(*)$ forces to set $q=1$
in the group $G_D$. Thus $\Z[G_D]$ cannot be a subring of $\Z[G]$
and so a naive attempt to obtain a representation of the braid group
$\B_{0,n}(\Sigma)$ using a covering of $B_{n,k}(\Sigma)$
corresponding to any epimorphism $\phi :\B_{n,k}(\Sigma)\to G$ cannot
give an extension of any homology linear representation of the
classical braid groups.

\section{A family of proposed representations}
As we have seen in the previous section, we are forced to take a
rather small covering of $B_{n,k}(\Sigma)$ in order that the
condition $(*)$ is satisfied, that is, the braid action commutes with
covering transformations so that it preserves the coefficient. A
remedy that we propose in this article is to use the same
configuration space $B_{n,k}(\Sigma)$ with an extended coefficient
ring so that we have some room to adjust coefficients to make the
braid action compatible with the coefficients. This remedy is a
reasonable thing to do if we hope to construct an extension of
homology linear representations of the classical braid groups.
Indeed, we successfully obtain an extension that seems the most
general among ones obtained from coverings of $B_{n,k}(\Sigma)$.

\subsection{Existence of extension of homology linear representation}
We first consider the intertwining braid group $\B_{n;k}(\Sigma)$. Note that $\B_{n;k}(\Sigma)$ is a candidate of group extension of $\B_{n,k}(\Sigma)$ by Lemma~\ref{lem:splits},
and $\B_{0,n}(\Sigma)$ acts on $\B_{n;k}(\Sigma)$ by right multiplication and acts on $\B_{n,k}(\Sigma)$ by conjugate because $\B_{n;k}(\Sigma)$ is the semidirect product of $\B_{0,n}(\Sigma)$ and $\B_{n,k}(\Sigma)$.

Let $H_\Sigma$ be the abstract group depending only on $k$ which admits the following presentation: for $k\ge 2$,
\begin{itemize}
\item Generators : $q, t, \bar m_1,\dots,\bar m_g, \bar\ell_1,\dots,\bar\ell_g,
m_1,\dots,m_g,\ell_1,\dots,\ell_g$
\item Relations : All generators commute except for the following
$$
[m_r,\ell_r]=t^2, \quad
[\bar m_r,\ell_r]=[m_r,\bar \ell_r]=q
$$
\end{itemize}

Let $\psi_\Sigma:\B_{n;k}(\Sigma)\to H_\Sigma$ be the epimorphism onto a group $H_\Sigma$ defined by
$$
\psi_\Sigma(\sigma_i)=t,\quad \psi_\Sigma(\zeta_j)=q,\quad \psi_\Sigma(\bar\sigma_m)=1
$$
and
$$
\psi_\Sigma(\mu_r)=m_r,\quad \psi_\Sigma(\lambda_r)=\ell_r,\quad
\psi_\Sigma(\bar\mu_r)=\bar m_r,\quad \psi_\Sigma(\bar\lambda_r)=\bar\ell_r
$$
where $1\le i\le k-1$, $1\le j\le n, 1\le m\le n-1$ and $1\le r\le g$.
If $k=1$, then we define $H_\Sigma$ to be the quotient of the above group by $t=1$.
Then $H_D$ is isomorphic to $G_D$ defined earlier for all $k\ge 1$, and is a subgroup of $H_\Sigma$ for any $\Sigma$ and $k\ge 1$. Even though $H_\Sigma$ (or $H_D$) depends on whether $k=1$ or $k\ge 2$, our notation does not indicate it for the sake of simplicity.

Let $\phi_\Sigma:\B_{n,k}(\Sigma)\to G_\Sigma$ be the restriction of $\psi_\Sigma$ to $\B_{n,k}(\Sigma)$ onto $G_\Sigma$ that denotes the normal subgroup $\psi_\Sigma(\B_{n,k}(\Sigma))$ of $H_\Sigma$ generated by $\{q, t, m_1,\dots,m_g,\ell_1,\dots,\ell_g\}$.
Then we can find the covering $p:\tilde{B}_{n,k}(\Sigma)\to B_{n,k}(\Sigma)$ corresponding to $\ker \phi_\Sigma$. Since the braid group $\B_{0,n}(\Sigma)$ embeds into the mapping class group of punctured surface $\Sigma_n$, a braid $\beta\in\B_{0,n}(\Sigma)$ determines a homeomorphism $\bar\beta:B_{n,k}(\Sigma)\to B_{n,k}(\Sigma)$.
Recall that the induced homomorphism $\bar\beta_*$ on $\B_{n,k}(\Sigma)$ is in fact the same as the conjugation by $\iota(\beta)$ where $\iota:\B_{0,n}(\Sigma)\to \B_{n;k}(\Sigma)$ is the splitting map in Lemma~\ref{lem:splits}.

\begin{lem}\label{lem:surfaceliftingcriteria}
Using the notations above,
the homeomorphism $\bar\beta:B_{n,k}(\Sigma)\to B_{n,k}(\Sigma)$ lifts to a homeomorphism $\tilde\beta:\tilde B_{n,k}(\Sigma)\to\tilde B_{n,k}(\Sigma)$ for any $\beta\in\B_{0,n}(\Sigma)$ and the restriction  $\phi_{\Sigma}$ of $\psi_{\Sigma}$ satisfies $\beta_\sharp\phi_{\Sigma}=\phi_{\Sigma}\bar\beta_*$
\end{lem}

\begin{proof}
By the lifting criteria, $\bar\beta$ lifts to $\tilde\beta$ if and only if $\bar\beta_*(\ker\phi_\Sigma)\subset \ker\phi_\Sigma$ if and only if there is an induced automorphism $\beta_\sharp$ on $G_\Sigma$ given by $\beta_\sharp\phi_\Sigma=\phi_\Sigma\bar\beta_*$.
Thus it suffices to show that $\phi_{\Sigma}\bar\beta_*(W)=1$ for any $W\in\ker \phi_{\Sigma}$
and $\beta\in\B_{0,n}(\Sigma)$. Let $W$ be a word in generators $\{\mu_i,\lambda_i,\sigma_i,\zeta_i\}$ of $\B_{n,k}(\Sigma)$.
Since the presentation for $H_\Sigma$ shows that any two elements are commutative up to multiplications by central elements $q$ and $t$, we have
$$
\phi_\Sigma\left(W\right)=W(\mu_i\leftarrow m_i,\lambda_i\leftarrow\ell_i,\sigma_i\leftarrow t,\zeta_i\leftarrow q)=q^c t^d
\prod m_i^{a_i} \ell_i^{b_i}
$$
where 
$W(\{x_i\leftarrow y_i\})$ denote the word obtained from $W$ by replacing the generators $x_i$'s by $y_i$'s.

Suppose $\phi_{\Sigma}(W)=1$. Then $a_i=b_i=0$ for all $1\le i\le
g$. Thus for generators $\sigma_r, \mu_r, \lambda_r$ for $\B_{0,n}(\Sigma)$, we have
\begin{align*}
    \phi_{\Sigma}((\bar\sigma_r)_*\left( W\right))
    &=\phi_{\Sigma}(W(\zeta_r\leftarrow \zeta_r \zeta_{r+1}\zeta_r^{-1}, \zeta_{r+1}\leftarrow\zeta_r))\\
&=W(\mu_i\leftarrow m_i,\lambda_i\leftarrow\ell_i,\sigma_i\leftarrow t,\zeta_i\leftarrow q)=1,\\
    \phi_{\Sigma}((\bar \mu_r)_*\left( W\right))
&=\phi_\Sigma(W(\lambda_r\leftarrow\lambda_r\mu_r\zeta_1^{-1}\mu_r^{-1}))\\
&=W(\mu_i\leftarrow m_i,\lambda_i\leftarrow\ell_i,\lambda_r\leftarrow q^{-1}\ell_r,\sigma_i\leftarrow t,\zeta_i\leftarrow q)\\
&=q^{-b_r}W(\mu_i\leftarrow m_i,\lambda_i\leftarrow\ell_i,\sigma_i\leftarrow t,\zeta_i\leftarrow q)=1,\\
    \text{and}\quad\phi_{\Sigma}((\bar \lambda_r)_*\left( W\right))
&=\phi_\Sigma(W(\mu_r\leftarrow \lambda_r\zeta_1\lambda_r^{-1}\mu_r\zeta_1\lambda_r\zeta_1^{-1}\lambda_r^{-1}))\\
&=W(\mu_i\leftarrow m_i,\mu_r\leftarrow qm_r, \lambda_i\leftarrow\ell_i,\sigma_i\leftarrow t,\zeta_i\leftarrow q)\\
&=q^{a_i}W(\mu_i\leftarrow m_i,\lambda_i\leftarrow\ell_i,\sigma_i\leftarrow t,\zeta_i\leftarrow q)=1.
\end{align*}
Therefore
$\phi_\Sigma(\bar\beta_*(W))=1$ and so $\beta_\sharp(\phi_{\Sigma}(\alpha))=\phi_{\Sigma}(\bar\beta_*(\alpha))$ holds for all $\alpha\in\B_{n,k}(\Sigma)$.
\end{proof}

By the above lemma, we now have a $\Z$-module automorphism $\tilde\beta_*$ on $H_k^{BM}(\tilde B_{n,k}(\Sigma))$.
Notice that $\tilde\beta_*$ is not necessarily a $\Z[G_\Sigma]$-module
homomorphism since the condition $(*)$ in \S2.1 may not hold, that is, the automorphism $\beta_\sharp$ of $G_\Sigma$ may not be the identity.

On the other hand, $\B_{0,n}(\Sigma)$ acts on $\B_{n;k}(\Sigma)$ by the right multiplication and so
there is an induced action of $\beta$ on $H_\Sigma$ given by $\beta\cdot h=h\psi_\Sigma(\beta)$ for $h\in H_\Sigma$.
It is possible to alter the induced action by multiplying with a certain function $\chi$ from $\B_{0,n}(\Sigma)$ to the centralizer of $G_\Sigma$ in $H_\Sigma$. We will closely discuss this possibility in Thoerem~\ref{thm:character}.

Using the $\Z$-module automorphism $\tilde\beta_*$ and the action on $\B_{n;k}(\Sigma)$ by $\B_{0,n}(\Sigma)$, we construct a $\Z[H_\Sigma]$-module automorphism $\beta\otimes\tilde\beta_*$ on $\Z[H_\Sigma]\otimes_{\Z[G_\Sigma]}H_k^{BM}(\tilde B_{n,k}(\Sigma))$ by
$$
(\beta\otimes\tilde\beta_*)(h\otimes c)=(\beta\cdot h)\otimes\tilde\beta_*(c)
$$
for $h\in H_\Sigma$ and $c\in H_k^{BM}(\tilde B_{n,k}(\Sigma))$.

\begin{thm}
\label{thm:representation}
Let $\Sigma$ be a compact, connected, oriented 2-dimensional manifold
with non-empty boundary.
Let $H_\Sigma$ be the group (depending on $k$) and $\psi_\Sigma:\B_{n;k}(\Sigma)\to H_\Sigma$ be the epimorphism given in
the discussion above and let $\phi_\Sigma:\B_{n,k}(\Sigma)\to G_\Sigma=\phi_\Sigma(\B_{n,k}(\Sigma))$ be the restriction of $\psi_\Sigma$.
Then there is a homomorphism
$$
\Phi_k:\B_{0,n}(\Sigma)\to
\aut_{\Z[H_\Sigma]}\left(\Z[H_\Sigma]\otimes_{\Z[G_\Sigma]}H_k^{BM}(\tilde B_{n,k}(\Sigma))\right)
$$
defined by
$$\Phi_k(\beta)=\beta\otimes\tilde\beta_*$$
where the action of $\beta$ on $H_\Sigma$ is given by
$\beta\cdot h=h\psi(\beta)$ for $h\in H_\Sigma$.

Moreover, this family $\Phi_k$ of representation is an extension of homology linear representation of the classical braid group $\B_{0,n}(D)$ in the sense of Definition~\ref{defn:extension}.
\end{thm}
\begin{proof}
Clearly $\Phi_k$ is a group homomorphism.
To see the well-definedness and the $\Z[H_\Sigma]$-linearity of $\Phi_k(\beta)$,
we assert that
$$\beta\cdot(hh')=(\beta\cdot h)\beta_\sharp(h').$$
for all $h\in H_\Sigma, h'\in G_\Sigma$.
Then
\begin{align*}
(\beta\otimes\tilde\beta_*)(h\otimes h'c)&=(\beta\cdot h)\otimes\tilde\beta_*(h'c)=(\beta\cdot h)\otimes \beta_\sharp(h')\tilde\beta_*(c)\\
&=(\beta\cdot h)\beta_\sharp(h')\otimes\tilde\beta_*(c)=
(\beta\otimes\tilde\beta_*)(hh'\otimes c)=hh'(\beta\otimes\tilde\beta_*)(1\otimes c)
\end{align*}
for $c\in H^{BM}_k(\tilde B_{n,k}(\Sigma))$.
Here the last equality is clear by the definition of the action by $\beta$.

To show the assertion, choose $\alpha\in\B_{n,k}(\Sigma)$ such that $\phi_\Sigma(\alpha)=h'$.
By Lemma~\ref{lem:surfaceliftingcriteria}, we have
$$\beta_\sharp(\phi_\Sigma(\alpha))=\phi_\Sigma(\bar\beta_*(\alpha))
=\psi_\Sigma(\bar\beta_*(\alpha))
=\psi_\Sigma(\beta^{-1}\alpha\beta)
=\psi_\Sigma(\beta)^{-1}\phi_\Sigma(\alpha)\psi_\Sigma(\beta)$$
Thus
$$\beta\cdot(hh')=hh'\psi_\Sigma(\beta)
=(h\psi_\Sigma(\beta))(\psi_\Sigma(\beta)^{-1}h'\psi_\Sigma(\beta))
=(\beta\cdot h)\beta_\sharp(h').$$

To show that $\Phi_k$ is
an extension of homology linear representations of the classical
braid groups, we regard an $n$ punctured disk $D_n$ as a subspace of
$\Sigma_n$. Then the configuration space $B_{n,k}(D)$ is a
subspace of the configuration space $B_{n,k}(\Sigma)$. For the
covering $p:\tilde B_{n,k}(\Sigma)\to B_{n,k}(\Sigma)$ corresponding
to $\phi_\Sigma$, a connected component of
$p^{-1}(B_{n,k})$ is denoted by $\tilde B_{n,k}(D)$. Since $G_D$ embeds into $G_\Sigma$, $p|_{\tilde
B_{n,k}(D)}:\tilde B_{n,k}(D)\to B_{n,k}(D)$ is the covering over $B_{n,k}(D)$
corresponding to $\psi|_{\B_{n,k}(D)}:\B_{n,k}(D)\to G_D$.
In fact, $H_k^{BM}(\tilde
B_{n,k}(D))$ is a submodule of $H_k^{BM}(\tilde B_{n,k}(\Sigma))$ as
$\Z[G_D]$-modules. One can see this more explicitly in the proof of
Lemma~\ref{lem:dimension}.
Each braid $\beta\in\B_{0,n}(D)$ gives a
$\Z[H_D]$-module automorphism $\beta\otimes\tilde\beta_*$ on
$\Z[H_D]\otimes_{\Z[G_D]} H_k^{BM}(\tilde B_{n,k}(D))$. Since
$G_D=H_D$, this automorphism is the same as $\tilde\beta_*$ on
$H_k^{BM} (\tilde B_{n,k}(D))$, which is the homology linear
representation of the classical braid group.
\end{proof}
In \S4, we will show that if one wants to obtain the result similar to the above theorem, the extension $\B_{n;k}(\Sigma)$ of $\B_{n,k}(\Sigma)$ is determined
uniquely up to redundant coefficient extension, and the quotient $H_\Sigma$ is uniquely determined for $k\ge 3$ and
is the simplest for $k\ge 1$ in the sense that any proper quotient of $H_\Sigma$ does not contain $G_D$ properly.

\subsection{Computation of proposed representations}
We now compute an explicit matrix forms of the representations
described in Theorem~\ref{thm:representation}, which turn out
to be extensions of the Burau representation and Lawrence-Krammer-Bigelow
representation of the classical braid groups. The following lemma
and its proof show not only that $H_k^{BM}(\tilde B_{n,k}(\Sigma))$
is a free $\Z[G_{\Sigma}]$-module but also how to choose a basis.  The lemma
is an extension of the corresponding lemma on a disk by Bigelow
\cite{Big4} and we borrow the main idea of his proof.

\begin{lem}\label{lem:dimension}
The homology group $H_\ell^{BM}(\tilde B_{n,k}(\Sigma))$ is the direct sum of
$$
\binom{2g+n+k-2}{k}$$ copies of $\Z[G_{\Sigma}]$ for $\ell=k$ and trivial otherwise.
\end{lem}
\begin{proof}
Let $d$ be a metric on $\Sigma$ that can be either hyperbolic or
Euclidean. Suppose punctures $z_1^0,\ldots,z_n^0$ lie on a geodesic.
And let $\gamma_j$ be the geodesic segment joining $z_j^0$ and
$z_{j+1}^0$ for $1\le j\le n-1$. For $1\le i\le g$, let $\alpha_i$,
$\beta_i$ be geodesic loops based at $z_1^0$ that represent the
meridian and the longitude of the $i$-th handle so that
$\alpha_i$'s, $\beta_i$'s, and $\gamma_j$'s are mutually disjoint.
Let $\Gamma$ be the union of all of these arcs so that $\Gamma_n$ is
a disjoint union of open $2g+n-1$ geodesic segments. Consider
$$B_\Gamma=B_{n,k}(\Gamma)=
\left\{\{z_1,\ldots,z_k\}\subset \Gamma_n\right\}.$$ Then it is not
hard to see $B_\Gamma$ is homeomorphic to a disjoint union of
$\binom{2g+n+k-2}{k}$ open $k$-balls that can be parameterized by
$(2g+n-1)$-tuples $(r_1,\ldots,r_{2g+n-1})$ of nonnegative integers
that add up to $k$ so that the $i$-th segment of $\Gamma_n$ contains
$r_i$ points from $\{z_1,\ldots,z_k\}$.

Let $p:\tilde B_{n,k}(\Sigma)\to B_{n,k}(\Sigma)$ be the covering
corresponding to the epimorphism $\phi_{\Sigma}:B_{n,k}(\Sigma)\to
G_{\Sigma}$. We will be done if we show that
$$H_\ell^{BM}(p^{-1}(B_\Gamma))\to H_\ell^{BM}(\tilde
B_{n,k}(\Sigma))$$ induced by the inclusion is an isomorphism since
$H_\ell^{BM}(B_\Gamma)$ is isomorphic to the direct sum of
$\binom{2g+n+k-2}{k}$ copies of $H_\ell(D^k, S^{k-1})$.

Define the family of compact subsets $A_\epsilon$ of $\Sigma_n$
defined by
$$A_\epsilon=\left\{\{z_1,\dots,z_k\}\in B_{n,k}(\Sigma)\mid
d(z_i, z_j)\ge \epsilon \mbox{ for }i\neq j,\: d(z_i, z_j^0)\ge
\epsilon \mbox{ for all }i,j\right\}.$$
 Since any compact subset of
$B_{n,k}(\Sigma)$ is contained in $A_\epsilon$ for sufficiently
small $\epsilon>0$, it suffices to show that
$$H_\ell(p^{-1}(B_\Gamma),p^{-1}(B_\Gamma - A_\epsilon))
\to H_\ell(\tilde B_{n,k}(\Sigma),p^{-1}(B_{n,k}(\Sigma)-
A_\epsilon))$$
is an isomorphism.

Let $\Sigma_\epsilon \subset \Sigma$ be the closed $\epsilon$
neighborhood of $\Gamma$ and let
$B_\epsilon=B_{n,k}(\Sigma_\epsilon)$. Then the obvious homotopy
collapsing from $B_{n,k}(\Sigma)$ to $B_\epsilon$ gives the
isomorphism
$$H_\ell(p^{-1}(B_\epsilon),p^{-1}(B_\epsilon- A_\epsilon))
\to
H_\ell(\tilde B_{n,k}(\Sigma),p^{-1}(B_{n,k}(\Sigma)- A_\epsilon)).
$$
Let $B$ be the set of $\{x_1,\dots,x_k\}\in B_\epsilon$ such that
for each $x_i$ there exists the unique nearest point in $\Gamma_n$.
Then $B$ is open and contains $A_\epsilon\cap B_\epsilon$. By
excision, the inclusion induces an isomorphism
$$H_\ell(p^{-1}(B), p^{-1}(B- A_\epsilon))\to
H_\ell(p^{-1}(B_\epsilon),p^{-1}(B_\epsilon- A_\epsilon)).$$ Finally,
the obvious deformation retract from $B$ to $B_\Gamma$ gives an
isomorphism
$$
H_\ell(p^{-1}(B), p^{-1}(B- A_\epsilon))\to H_\ell(p^{-1}(B_\Gamma),
p^{-1}(B_\Gamma- A_\epsilon)).
$$
that completes the proof.
\end{proof}

We remark that $B_\epsilon=B_{n,k}(\Sigma_\epsilon)$ and
$B_\Gamma=B_{n,k}(\Gamma)$ do not have the same homotopy type even
though $\Gamma$ is a deformation retract of $\Sigma_\epsilon$. This
is because $\Gamma$ is a 1-dimensional complex and movements of
points in $\Gamma$ avoiding collision are more restricted.

Let $I(n,k,g)$ be the set of $(2g+n-1)$-tuples of nonnegative
integers that add up to $k$. The proof of Lemma~\ref{lem:dimension}
shows that a typical basis for $H_k^{BM}(\tilde B_{n,k}(\Sigma))$ as
a free $\Z[G_{\Sigma}]$-module can be indexed by the set $I(n,k,g)$. The
proof also shows that a basis for $H_k^{BM}(\tilde B_{n,k}(D))$ can be
chosen as a subset of a basis for $H_k^{BM}(\tilde
B_{n,k}(\Sigma))$. Thus the homology linear representations for
$\B_{0,n}(D)$ appears in matrix forms of proposed representations for
$\B_{0,n}(\Sigma)$ as minors.

\begin{figure}[ht]
\begin{center}
{\subfigure[Fork $F$]{\makebox[6cm]{\psfrag{A}{$\beta_1$}\psfrag{B}{$\alpha_2^2$}\psfrag{C}{$\gamma_2^3$}\includegraphics[scale=0.6]{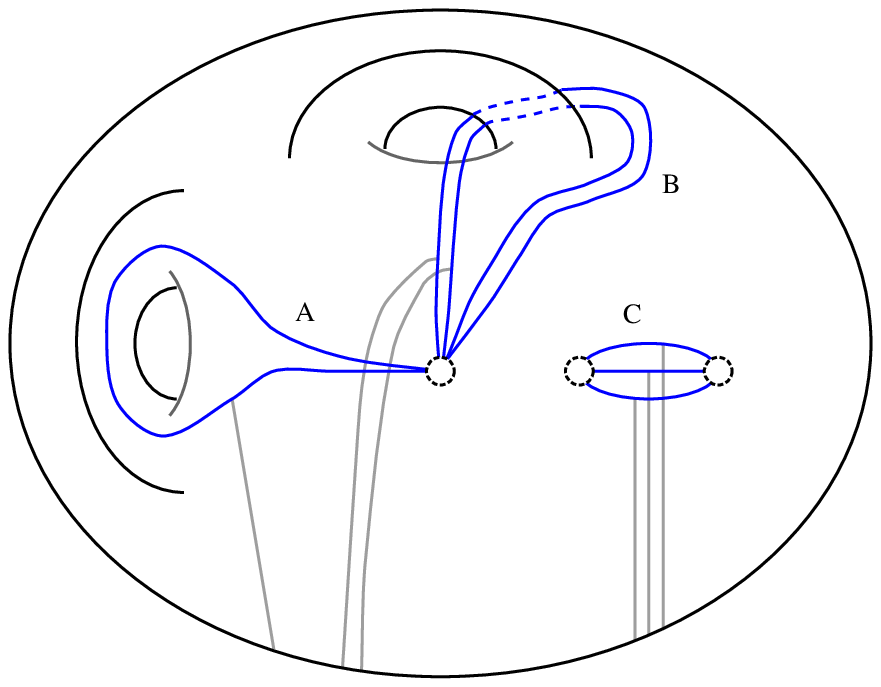}}}
\subfigure[Dual noodle $N$]{\makebox[6cm]{\psfrag{A}{$\beta_1^*$}\psfrag{B}{$(\alpha_2^2)^*$}\psfrag{C}{$(\gamma_2^3)^*$}\includegraphics[scale=0.6]{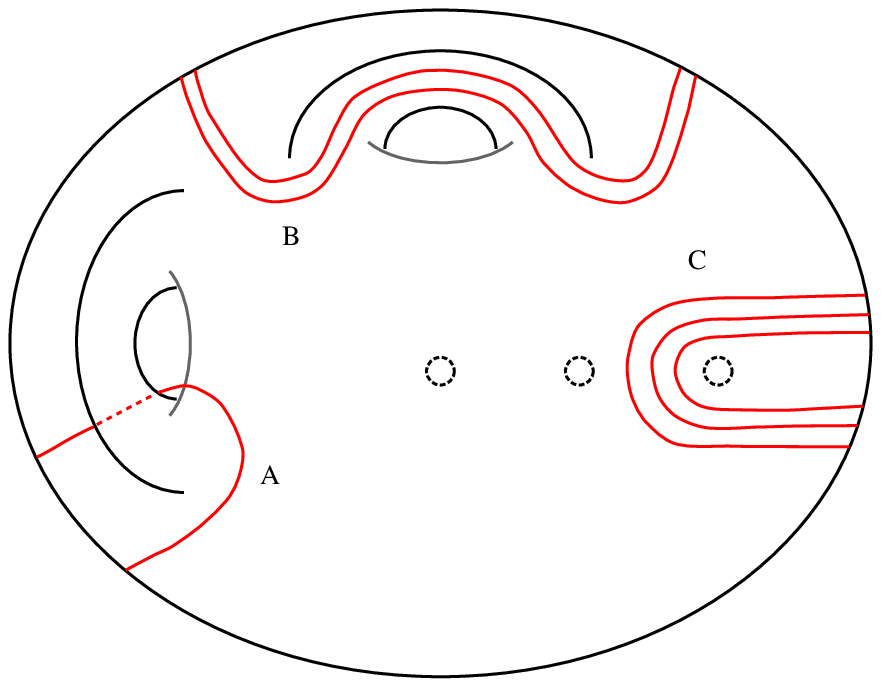}}}
}
\end{center}
\caption{An example of a fork and its dual noodle}
\label{fig:fork-noodle}
\end{figure}

By the work of Krammer \cite{Kra1,Kra2} and Bigelow
\cite{Big2,Big3,Big4}, it has been known that there is a natural and
useful way of describing a basis geometrically. Recall the loops
$\alpha_i, \beta_i$ and the arcs $\gamma_j$ from the proof of
Lemma~\ref{lem:dimension}. For $(r_1,\ldots,r_{2g+n-1})\in
I(n,k,g)$, choose $r_i$ disjoint duplicates of $\alpha_i$ or
$\beta_{i-g}$ or $\gamma_{i-2g}$ if $1\le i\le g$ or $g+1\le i\le
2g$ or $2g+1\le i\le 2g+n-1$, respectively. For each $i$, join these
$r_i$ disjoint duplicates to $\partial\Sigma$ by mutually  disjoint
arcs (that determine a basing). This geometric object is called a
{\em fork}. In fact, a fork uniquely determines a $k$-cycle in
$H_k^{BM}(\tilde B_{n,k}(\Sigma))$ by lifting the Cartesian product
of $k$ curves together with basing arcs in the fork. Notice that the
basing is required to have a unique lift. For example, the fork
corresponding to $(0,2,1,0,0,3)\in I(3,6,2)$ looks like the set of
curves on the left of Figure~\ref{fig:fork-noodle}.

As Bigelow showed for the case of the disk in \cite{Big4}, the
Poincar\'e duality, the universal coefficient theorem, and
Lemma~\ref{lem:dimension} imply that the ordinary relative homology
$H_k(\tilde B_{n,k}(\Sigma),\partial\tilde B_{n,k}(\Sigma))$ is the
dual space of the Borel-Moore homology $H_k^{BM}(\tilde
B_{n,k}(\Sigma))$ in the sense that there is a
nonsingular sesquilinear pairing
$$\langle\: , \rangle:H_k^{BM}(\tilde B_{n,k}(\Sigma))\times
H_k(\tilde B_{n,k}(\Sigma), \partial\tilde B_{n,k}(\Sigma))\to \Z[S]$$
where $S$ is a skew field containing $\Z[G_{\Sigma}]$. In fact, the
group $G_{\Sigma}$ is bi-ordered and so it can embed into a skew
field such as the Malcev-Neumann power series ring \cite{Mal, Neu}.
Explicitly, the above pairing is defined as follows:
for cycles $F\in H_k^{BM}(\tilde B_{n,k}(\Sigma)), N\in H_k(\tilde B_{n,k}(\Sigma),\partial\tilde B_{n,k}(\Sigma))$ in each homology groups
\begin{displaymath}
\langle F, N
\rangle=\sum_{y\in G_\Sigma} y(F, y N)
\end{displaymath}
where $(\: , )$ counts the intersection number.

Let $\alpha_1^*,\ldots,\alpha_g^*$, $\beta_1^*,\ldots,\beta_g^*$,
$\gamma_1^*,\ldots,\gamma_g^*$ be pairwise disjoint arcs which start
and end at $\partial\Sigma$ and $\alpha_i^*$ (or $\beta_i^*$, or
$\gamma_j^*$) intersects only $\alpha_i$ (or $\beta_i$, or
$\gamma_j$) once transversely. A basis of $H_k(\tilde
B_{n,k}(\Sigma),\partial\tilde B_{n,k}(\Sigma))$ by duplicating
$\alpha_i^*$'s or $\beta_i^*$'s or $\gamma_j^*$'s depending on a
given $(2g+n-1)$-tuples in $I(n,k,g)$. This geometric object is
called a {\em noodle}. In fact, a noodle uniquely determines a
relative $k$-cycle in $H_k(\tilde B_{n,k}(\Sigma),\partial\tilde
B_{n,k}(\Sigma))$ by lifting the Cartesian product of $k$ arcs in
the noodle. For example, the noodle corresponding to
$(0,2,1,0,0,3)\in I(3,6,2)$ looks like the set of arcs on the right
of Figure~\ref{fig:fork-noodle}.

For  a $k$-cycle $F$ determined by a fork, and a relative $k$-cycle
$N$ determined by a noodle, the sesquilinear pairing $\langle F,
N\rangle$ computes algebraic intersections between them with $\Z[G_{\Sigma}]$
coefficients. This pairing can easily be computed by recording
intersections between the fork and the noodle on $\Sigma$. The basis
determined by folks and the basis determined by noodles are dual
with respect to the pairing.

In the case of a disc, Bigelow \cite{Big4} showed that this pairing
is invariant under the action by $\B_{0,n}(D)$. But in the case of a
surface $\Sigma$ of genus $\ge 1$, it can not be invariant under the
action by $\B_{0,n}(\Sigma)$. In fact, the pairing cannot be
preserved by any braid group action given by a representation $\Psi$
into $\aut_{\Z[G_\Sigma]}\left(H_k^{BM}(\tilde B_{n,k}\right)$.
Suppose it is preserved, that is, for any
$\beta\in\B_{0,n}(\Sigma)$, a $k$-cycle $F$ determined by a fork,
and a relative $k$-cycle $N$ determined by a noodle,
$$
    \langle F, N \rangle = \langle \Psi(\beta)(F),
    \Psi(\beta)(N)\rangle,
$$
then for any $y\in G_{\Sigma}$,
\begin{align*}
    y\langle F, N\rangle &=\langle y F, N \rangle \\
    &= \langle \Psi(\beta)(y F), \Psi(\beta)(N) \rangle\\
    &=\beta_\sharp(y)\langle \Psi(\beta)(F), \Psi(\beta)(N) \rangle\\
    &=\beta_\sharp(y)\langle F, N \rangle.
\end{align*}
The property $y=\beta_\sharp(y)$ for all $y\in G_\Sigma$ forces to set $q=1$ in $G_{\Sigma}$
and so it was abandoned.

Nonetheless, we can extend this pairing to
$$
    \langle\: , \rangle_{H_{\Sigma}}:\Z[H_{\Sigma}]
    \otimes_{\Z[G_{\Sigma}]} H_k^{BM}(\tilde B_{n,k}(\Sigma))
    \times\Z[H_{\Sigma}]\otimes_{\Z[G_{\Sigma}]}H_k(\tilde B_{n,k}(\Sigma),
    \partial \tilde B_{n,k}(\Sigma))\to S'
$$
defined by
$$
    \left\langle \sum_i g_i F_i, \sum_j h_j N_j \right\rangle_{H_{\Sigma}}
    = \sum_{i,j} g_i \langle F_i, N_j\rangle h_j^{-1}
$$
where $g_i, h_j\in \Z[H_{\Sigma}]$ and $S'$ is the skew field
containing $\Z[H_{\Sigma}]$. Note that this extended pairing cannot
be invariant under the braid group action given by $\Phi_k$ either.
But it can be used to compute proposed representations $\Phi_k$
explicitly. The following theorem summarizes the above discussion.

\begin{thm}\label{thm:computation}
Let $F_i$'s and $N_i$'s be $k$-cycles and relative
$k$-cycles in dual bases determined by forks and noodles.
Then for
each $\beta\in\B_{0,n}(\Sigma)$, $\Phi_k(\beta)$ is represented by a
matrix with respect to the basis $\{F_i\mid 1\le
i\le\binom{2g+n+k-2}{k}\}$ whose $(i,j)$-th entry is given by
$$\psi_\Sigma(\beta)\left\langle \tilde{\beta}(F_i), N_j
\right\rangle_{H_{\Sigma}}$$ which is an element of $\Z[H_{\Sigma}]$
rather than of $S'$.
\end{thm}

 As an example, we will show the
matrix form of the representation $\Phi_1$ of the 3-braid group
$\B_{0,3}(\Sigma)$ which is an extension of the Burau representation
where $\Sigma=\Sigma(2,1)$. Since $k=1$, the basis of
$H_1^{BM}(\tilde B_{3,1}(\Sigma))$ determined by forks can be
expressed by $\{\gamma_1, \gamma_2, \alpha_1, \alpha_2, \beta_1,
\beta_2\}$ and similarly the dual basis of $H_1(\tilde
B_{3,1}(\Sigma),\partial\tilde B_{3,1}(\Sigma))$ determined by
noodles is written by $\{\gamma_1^*, \gamma_2^*,
\alpha_1^*, \alpha_2^*, \beta_1^*, \beta_2^*\}$.

\begin{figure}[ht]
\begin{center}
\subfigure[Fork $\beta_1$ and $\sigma_1$]{\psfrag{A}{$\beta_1$}\psfrag{C}{$\sigma_1$}\includegraphics[scale=0.6]{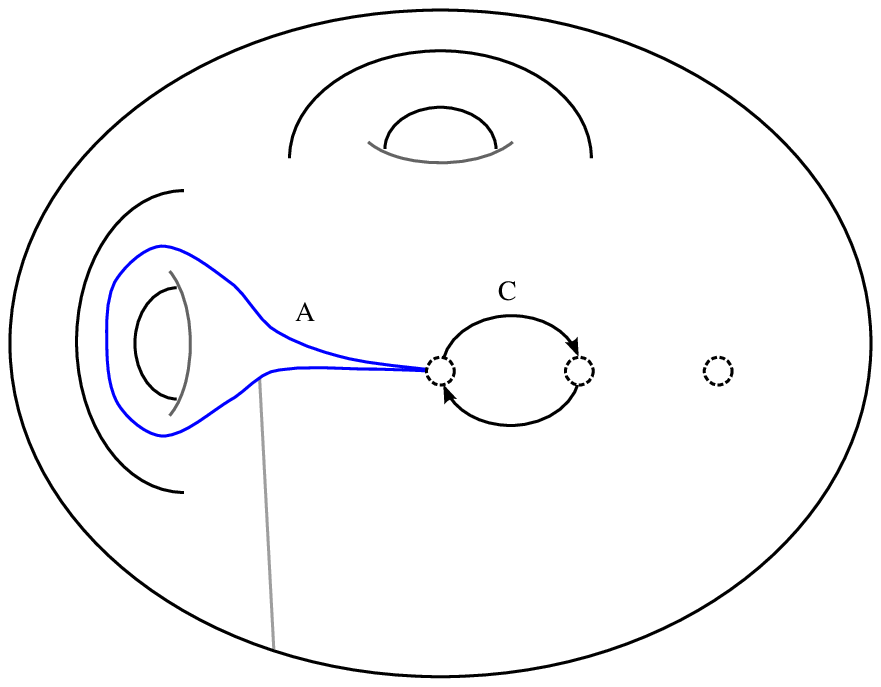}}
\subfigure[Fork and noodles]{\psfrag{A}{$(\tilde\sigma_1)_*(\beta_1)$}
\psfrag{B}{$\gamma_1^*$}
\psfrag{C}{$p_2$}
\psfrag{D}{$p_1$}
\psfrag{E}{$p_3$}
\psfrag{F}{$\beta_1^*$}
\includegraphics[scale=0.6]{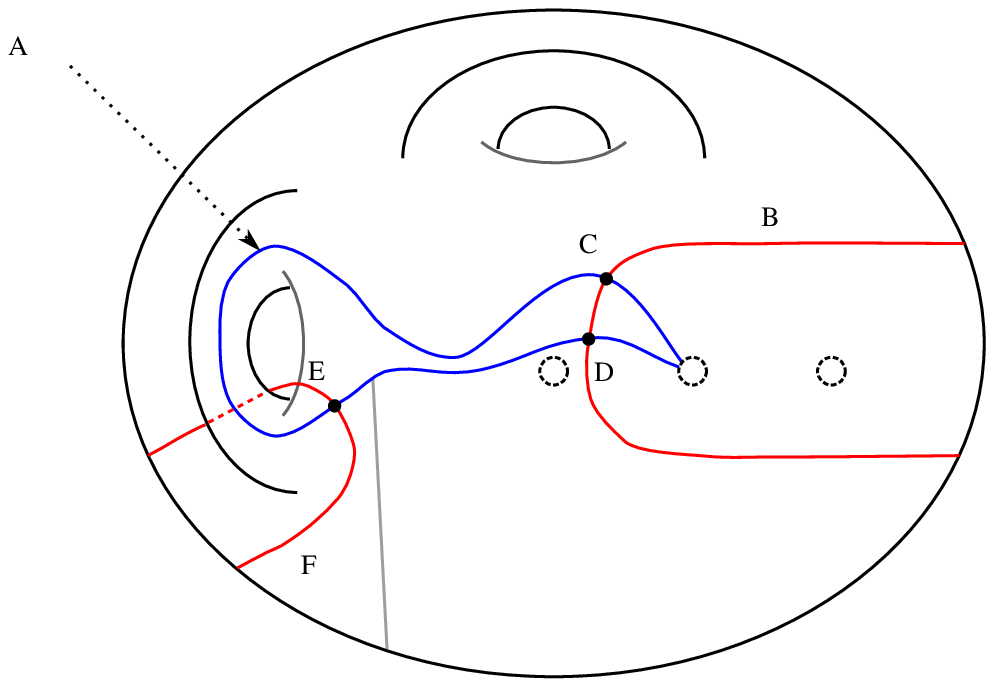}}
\end{center}
\caption{An example of the pairing of a fork and a noodle}
\label{fig:pairing}
\end{figure}

Figure~\ref{fig:pairing} shows the action of $\sigma_1$ on the fork
$\beta_1$. Also it shows intersection points $p_1, p_2$ with
$\gamma_1^*$ and $p_3$ with $\beta_1^*$. In the covering space, the
intersection point $p_1$ lies on the sheet transformed by $q$ since
the fork wraps a puncture, and $p_2$ lies on the sheet transformed
by $\ell_1 q$ since the fork contains the longitude of the first handle
and wraps a puncture. Besides we have the negative sign for $p_2$ since the orientation is switched. Finally $p_3$ lies on the sheet
containing the base point of the covering space. Therefore we have
$\Phi_1(\sigma_1)(\beta_1)=\beta_1 + q(1-\ell_1) \gamma_1$. By the
similar computation, we can obtain every entry of
$\Phi_1(\sigma_1)$.

\begin{center}
$$
\Phi_1(\sigma_1)=
\left(
\begin{array}{cc:cccc}
-q&1&q(1-m_1)&q(1-m_2)&q(1-\ell_1)&q(1-\ell_2)\\
0&1&0&0&0&0\\
\hdashline
\multicolumn{2}{c:}{\mathbf{0}}&\multicolumn{4}{c}{I_4}
\end{array}
\right)
$$
$$
\Phi_1(\sigma_2)=
\left(
\begin{array}{cc}
1&0\\
q&-q
\end{array}
\right)
\oplus I_4
$$
$$
\Phi_1(\mu_1)=
\bar m_1
\left(
\begin{array}{cc:cc:cc}
\multicolumn{2}{c:}{I_2}&\multicolumn{4}{c}{\mathbf{0}}\\
\hdashline
1&0&m_1q&q(m_2-1)&\ell_1-1&q(\ell_2-1)\\
0&0&0&1&0&0\\
\hdashline
\multicolumn{2}{c:}{\mathbf{0}}&\multicolumn{2}{c:}{\mathbf{0}}&\multicolumn{2}{c}{I_2}
\end{array}
\right)
$$
$$
\Phi_1(\mu_2)=
\bar m_2
\left(
\begin{array}{cc:cc:cc}
\multicolumn{2}{c:}{I_2}&\multicolumn{4}{c}{\mathbf{0}}\\
\hdashline
0&0&1&0&0&0\\
1&0&m_1-1&m_2q&\ell_1-1&\ell_2-1\\
\hdashline
\multicolumn{2}{c:}{\mathbf{0}}&\multicolumn{2}{c:}{\mathbf{0}}&\multicolumn{2}{c}{I_2}
\end{array}
\right)
$$
$$
\Phi_1(\lambda_1)=
\bar \ell_1
\left(
\begin{array}{cc:cc:cc}
\multicolumn{2}{c:}{I_2}&\multicolumn{4}{c}{\mathbf{0}}\\
\hdashline
\multicolumn{2}{c:}{\multirow{2}{*}{$\mathbf{0}$}}&q&0&\multicolumn{2}{c}{\multirow{2}{*}{$\mathbf{0}$}}\\
&&0&1\\
\hdashline
1&0&q(m_1q-1)&q(m_2-1)&\ell_1q&q(\ell_2-1)\\
0&0&0&0&0&1
\end{array}
\right)
$$
$$
\Phi_1(\lambda_2)=
\bar \ell_2
\left(
\begin{array}{cc:cc:cc}
\multicolumn{2}{c:}{I_2}&\multicolumn{4}{c}{\mathbf{0}}\\
\hdashline
\multicolumn{2}{c:}{\multirow{2}{*}{$\mathbf{0}$}}&1&0&\multicolumn{2}{c}{\multirow{2}{*}{$\mathbf{0}$}}\\
&&0&q\\
\hdashline
0&0&0&0&1&0\\
1&0&m_1-1&q(m_2q-1)&\ell_1-1&\ell_2q
\end{array}
\right)
$$
\end{center}

Similarly, we can compute the matrix form for $k=2$ which is the extension of Lawrence-Krammer-Bigelow representation. For $g=1, n=3$, by Lemma~\ref{lem:dimension}, we have $10\times 10$ matrix for each generator.
Fix a basis for $H^{BM}_2(\tilde B_{3,2}(\Sigma))$ as shown in Figure~\ref{fig:fork-noodle}.
Let $w_{1,1}=(0,0,2,0)$, $w_{1,2}=(0,0,1,1)$, $w_{2,2}=(0,0,0,2)$,
$a_{0,0}=(2,0,0,0)$, $a_{0,1}=(1,0,1,0)$, $a_{0,2}=(1,0,0,1)$,
$b_{0,0}=(0,2,0,0)$, $b_{0,1}=(0,1,1,0)$, $b_{0,2}=(0,1,0,1)$, and $z=(1,1,0,0)$ in $I(3,2,1)$.
Then the action of $\sigma_1$ on these basis is as follows:
\begin{align*}
\Phi_2(\sigma_1)(w_{1,1})&=tq^2 w_{1,1}\\
\Phi_2(\sigma_1)(w_{1,2})&=-tq w_{1,1}-q w_{1,2}\\
\Phi_2(\sigma_1)(w_{2,2})&=w_{1,1}+(1+t^{-1}) w_{1,2}+ w_{2,2}\\
\Phi_2(\sigma_1)(a_{0,0})&=a_{0,0}+q(1+t^{-1})(1-m_1 t)a_{0,1}
+q^2(m_1^2-(1+t)m_1+1)w_{1,1}\\
\Phi_2(\sigma_1)(a_{0,1})&=-q a_{0,1}+q^2t(m_1-1) w_{1,1}\\
\Phi_2(\sigma_1)(a_{0,2})&=a_{0,1}+a_{0,2}+qt(1-m_1) w_{1,1}+q(1-m_1) w_{1,2}\\
\Phi_2(\sigma_1)(b_{0,0})&=b_{0,0}+q(1+t^{-1})(1-\ell_1 t)b_{0,1}
+q^2(\ell_1^2-(1+t)\ell_1+1)w_{1,1}\\
\Phi_2(\sigma_1)(b_{0,1})&=-q b_{0,1}+q^2t(\ell_1-1) w_{1,1}\\
\Phi_2(\sigma_1)(b_{0,2})&=b_{0,1}+b_{0,2}+qt(1-\ell_1) w_{1,1}+q(1-\ell_1) w_{1,2}\\
\Phi_2(\sigma_1)(z)&=q(t^{-1}-t\ell_1)a_{0,1}+q(1-m_1)b_{0,1}
+q^2(1+m_1(\ell_1-1)-t \ell_1) w_{1,1}+z
\end{align*}
Note that the correspondence between the basis $\{v_{j,k}\}$ in \cite{Big2} and our basis is given as follows:
\begin{align*}
v_{1,2}&=-tq^{-4}w_{1,1}\\
v_{1,3}&=-tq^{-4}(w_{1,1}+q(1-t^{-1})w_{1,2}+q^2w_{2,2})\\
v_{2,3}&=-tq^{-2}w_{2,2}
\end{align*}
Then the action of $\Phi_2$ on the basis $\{v_{j,k}\}$ together with substitution $t\mapsto -t$ is exactly equal to that of Lawrence-Krammer-Bigelow representation in \cite{Big2}.

\section{Justification of proposed representation}
Beside the family of representations proposed in the previous section, we will investigate the possibility that there may be other representations of the surface braid groups that extends the homology linear representations of the classical braid groups. One may try to consider alternatives in the three ways---a group extension of $\B_{n,k}(\Sigma)$ other than $\B_{n;k}(\Sigma)$, a quotient group of $\B_{n;k}(\Sigma)$ other than $H_\Sigma$, and an action on $H_\Sigma$ by $\B_{0,n}(\Sigma)$ other than the right multiplication via the quotient map.

\subsection{Group extension of $\B_{n,k}(\Sigma)$}

In order to make an adjustment of coefficients in the most flexible
manner, we may try to find a group extension $\E_{n,k}(\Sigma)$ of
$\B_{n,k}(\Sigma)$ such that $\B_{0,n}(\Sigma)$ acts on
$\E_{n,k}(\Sigma)$ and it is as large as possible. If we regard
$\B_{0,n}(\Sigma)$ and $\B_{n,k}(\Sigma)$ as subgroups of some large
braid group $\B_{0,n+k+\ell}(\Sigma)$, then $\B_{0,n}(\Sigma)$ acts
naturally on $\B_{0,n+k+\ell}(\Sigma)$ as well as on
$\B_{n,k}(\Sigma)$ by conjugation. Thus we assume that
$\B_{n,k}(\Sigma)\subset\E_{n,k}(\Sigma)\subset\B_{0,n+k+\ell}(\Sigma)$
for some $\ell\ge 0$.

\begin{lem}
\label{lem:normalizer}
Let $\Sigma$ be a surface with non-empty boundary and
$\Sigma'$ be a collar neighborhood of $\partial\Sigma$. Let
$N(\B_{n,k}(\Sigma))$ denote the normalizer of $\B_{n,k}(\Sigma)$ in
$\B_{0,n+k+\ell}(\Sigma)$ for some $\ell\ge 0$. Then
$$N(\B_{n,k}(\Sigma))\cong\B_{n;k}(\Sigma)\times\B_{0,\ell}(\Sigma').$$
\end{lem}
\begin{proof} We first identify $\B_{n,k}(\Sigma)$ and $\B_{0,n}(\Sigma)$
to the corresponding subgroups of $\B_{0,n+k+\ell}(\Sigma)$ via the
embeddings which add trivial $\ell$ and $k+\ell$ strands,
respectively. Then we will show
$N(\B_{n,k}(\Sigma))=\B_{n;k}(\Sigma)\times\B_{0,\ell}(\Sigma')$ as
subgroups of $\B_{0,n+k+\ell}(\Sigma)$. It is clear that
$\B_{n;k}(\Sigma)\times\B_{0,\ell}(\Sigma')\subset
N(\B_{n,k}(\Sigma))$ since $\B_{n,k}(\Sigma)$ is a normal subgroup
of $\B_{n;k}(\Sigma)$ from the short exact sequence of
Lemma~\ref{lem:splits} and elements of $\B_{0,\ell}(\Sigma')$
commutes with those of $\B_{n,k}(\Sigma)$. Conversely let $\beta\in
N(\B_{n,k}(\Sigma))\subset\B_{0,n+k+\ell}(\Sigma)$. Any element
$\alpha\in\B_{n,k}(\Sigma)$ and its conjugate
$\beta^{-1}\alpha\beta\in\B_{n,k}(\Sigma)$ induce permutations that
preserve the sets $\left\{ 1,\dots,n \right\}$, $\left\{
n+1,\dots,n+k \right\}$ and $\left\{ n+k+1,\dots,n+k+\ell \right\}$.
It is easy to see that the induced permutation of $\beta$ itself
must fix the above three sets since $\alpha$ can be arbitrary in
$\B_{n,k}(\Sigma)$. Thus $\beta\in\B_{n+k;l}(\Sigma)$ and the split
exact sequence
$$
\xymatrix{
    1\ar[r]&\B_{n+k,l}(\Sigma)\ar[r]&\B_{n+k;l}(\Sigma)
    \ar[r]^{(\pi_{n+k})_*}&\B_{0,n+k}(\Sigma)\ar[r]&1
    }
$$
gives a unique decomposition $\beta=\beta_1\beta_2$ for $\beta_1\in
\B_{0,n+k}(\Sigma)$ and $\beta_2\in \B_{n+k,\ell}(\Sigma)$. In fact,
$\beta_1=(\pi_{n+k})_*(\beta)\in \B_{n;k}(\Sigma)$ since the epimorphism
$(\pi_{n+k})_*$ forgets the last $\ell$ strands or replaces them by
the trivial $\ell$-strand braid.

For any $\alpha\in \B_{n,k}(\Sigma)\subset \B_{0,n+k}(\Sigma)\subset
\B_{0,n+k+\ell}(\Sigma)$, we have
$(\pi_{n+k})_*(\beta_2^{-1}\alpha\beta_2)=\beta_2^{-1}\alpha\beta_2$
since $\beta_2^{-1}\alpha\beta_2\in \B_{0,n+k}$. On the other hand,
we have $(\pi_{n+k})_*(\beta_2^{-1}\alpha\beta_2)=\alpha$ since
$(\pi_{n+k})_*$ replaces the last $\ell$ strands by the trivial
braid. Thus we have $\beta_2^{-1}\alpha\beta_2=\alpha$. From the
presentation of $\B_{0,n+k+\ell}(\Sigma)$ in \S1.1, it is easy to
see that $\beta_2\in \B_{n+k,\ell}(\Sigma)$ must be a local braid in
order for $\beta_2$ to commute with every element of
$\B_{n,k}(\Sigma)$. Thus we have $\beta_2\in \B_{0,\ell}(\Sigma')$
where $\Sigma'$ is an annulus that is a collar neighborhood of
$\partial \Sigma$ in $\Sigma$. Consequently, we have shown
$N(\B_{n,k}(\Sigma))\subset\B_{n;k}(\Sigma)\times\B_{0,\ell}(\Sigma')$
\end{proof}

By the above lemma, the extension $\E_{n,k}(\Sigma)$ of
$\B_{n,k}(\Sigma)$ can be taken as a subgroup of
$\B_{n;k}(\Sigma)\times\B_{0,l}(\Sigma')$. We remark that
$\B_{n;k}(\Sigma)\times\B_{0,l}(\Sigma')$ is also a subgroup  of the
intertwining braid group $\B_{n;k+l}(\Sigma)$.

Then we follow the construction given in the discussion before Theorem~\ref{thm:representation} with $\E_{n,k}(\Sigma)$ replacing $\B_{n;k}(\Sigma)$.

Let $\psi:\E_{n,k}(\Sigma) \to H$ be an epimorphism onto a group $H$.
If we choose an action of $\B_{0,n}(\Sigma)$ on the extension $\E_{n,k}(\Sigma)$ then the action is carried over $H$ via $\psi$ and it is convenient to use the convention that $(\beta_1\beta_2)\cdot h=\beta_2\cdot(\beta_1\cdot h)$ for $\beta_1,\beta_2\in\B_{0,n}(\Sigma)$ and $h\in H$. In order to obtain a $\Z[H]$-module automorphism
$\beta\otimes\tilde\beta_*$ on $\Z[H]\otimes_{\Z[G]}H^{BM}_k(\tilde B_{n,k}(\Sigma))$ that is an extension of a homology linear representation  of the classical braid group, this induced action of $\B_{0,n}(\Sigma)$ on $H$ needs to satisfy the following two conditions:
\begin{itemize}
\item [(i)] Lifting criteria : $\beta_\sharp$ exists and $\beta_\sharp(\phi(\alpha))=\phi(\bar\beta_*(\alpha))$ for all $\alpha\in \B_{n,k}(\Sigma)$ where $\phi=\psi|_{\B_{n,k}(\Sigma)}$;
\item [(ii)] Linearity and compatibility :
$hh'(\beta\cdot 1)=\beta\cdot(hh')=(\beta\cdot h)\beta_\sharp(h')$
for all $h\in H, h'\in G=\phi(\B_{n,k}(\Sigma))$.
\end{itemize}
As in the proof of Theorem~\ref{thm:representation}, we then have
$$(\beta\otimes\tilde\beta_*)(h\otimes h'c)=
(\beta\otimes\tilde\beta_*)(hh'\otimes c)=
hh'(\beta\otimes\tilde\beta_*)(1\otimes c)$$
for all $h\in H$, $h'\in G$ and $c\in H_k^{BM}(\tilde B_{n,k}(\Sigma))$.

\begin{thm}\label{thm:uniqueness}
Suppose that there are an epimorphism $\psi:\E_{n,k}(\Sigma)\to H$ and an action of $\B_{0,n}(\Sigma)$ on $H$ satisfying the above two conditions. Let $\Psi_k$ be the representation obtained from $\psi$ and the action.
Then $$\Psi_k=1_{\Z[H]}\otimes_{\Z[H']}\Psi'_k$$
for a representation $\Psi'_k$ obtained from an epimorphism $\psi':\B_{n;k}(\Sigma)\to H'\subset H$ and an action $\B_{0,n}(\Sigma)$ on $H'$
where $1_{\Z[H]}$ is the identity map on $\Z[H]$.
\end{thm}
\begin{proof}
Let $H'=\left\{\beta\cdot1\in H|\beta\in\B_{0,n}(\Sigma)\right\}\phi(\B_{n,k}(\Sigma))$ and $\psi':\B_{n;k}(\Sigma)\to H'$ be a surjection defined by $\psi'(\beta)=\beta\cdot 1$ for $\beta\in\B_{0,n}(\Sigma)$ and $\psi'=\phi$ on $\B_{n,k}(\Sigma)$.
Then since
$$
\psi'(\beta_1\beta_2)=(\beta_1\beta_2)\cdot 1
=\beta_2\cdot(\beta_1\cdot 1)
=(\beta_1\cdot 1)(\beta_2\cdot 1)=\psi'(\beta_1)\psi'(\beta_2)
$$
for all $\beta_1,\beta_2\in\B_{0,n}(\Sigma)$, $\psi'$ becomes a homomorphism and
$\psi'$ preserves the semidirect product structure.
Also we have
$$
\phi'=\psi'|_{\B_{n,k}(\Sigma)}=\psi|_{\B_{n,k}(\Sigma)}=\phi
$$
and so $\phi$ and $\phi'$ induce the same homology group $H_k^{BM}(\tilde B_{n,k}(\Sigma))$ and two $\Z$-module automorphisms obtained from $\beta$ coincide.

Now consider two representations $\Psi_k$ and $\Psi'_k$ corresponding to
$\psi$ and $\psi'$, respectively. Then $\Psi_k(\beta)$ gives a
$\Z[H]$-homomorphism on $\Z[H]\otimes_{\Z[G]}H_k^{BM}(\tilde
B_{n,k}(\Sigma))$ and $\Psi'_k(\beta)$ gives a $\Z[H']$-homomorphism
on $\Z[H']\otimes_{\Z[G]}H_k^{BM}(\tilde B_{n,k}(\Sigma))$.
Since
$\Z[H]=\Z[H]\otimes_{\Z[H']}\Z[H']$, $\Psi_k(\beta)$ is a
$\Z[H]$-homomorphism on
$\Z[H]\otimes_{\Z[H']}\Z[H']\otimes_{\Z[G]}H_k^{BM}(\tilde
B_{n,k}(\Sigma))$ defined by
$$\Psi_k(\beta)(hh'\otimes c)=hh'(\beta\cdot1)\otimes\tilde\beta_*(c)=
h\otimes h'(\beta\cdot 1)\otimes\tilde\beta_*(c)$$ for all $h\in \Z[H],
h'\in \Z[H']$ and $c\in H_k^{BM}(\tilde B_{n,k}(\Sigma))$. This is equal to $1_{\Z[H]}\otimes_{\Z[H']}\Psi'_k(\beta)$.
\end{proof}

The above theorem implies that we may assume that $\B_{n;k}(\Sigma)\subset\E_{n,k}(\Sigma)$ without loss of generality.
Then
by Lemma~\ref{lem:normalizer}, $\E_{n,k}(\Sigma)=\B_{n;k}(\Sigma)\times \B$ for some subgroup $\B$ of $\B_{0,\ell}(\Sigma')$ and the above theorem says that any family of representations obtained by using $\E_{n,k}(\Sigma)$ is merely a trivial extension of the family of representations proposed in \S3.

\subsection{Quotient of $\B_{n;k}(\Sigma)$}
According to the scheme described in
Theorem~\ref{thm:representation}, it is important to find a
good epimorphism $\psi:\B_{n;k}(\Sigma)\to H$ onto some group $H$.

Since $\Sigma$ is not sphere, the inclusion $B_{n,k}(D)\hookrightarrow
B_{n,k}(\Sigma)$ induces a monomorphism $\B_{n,k}(D)\hookrightarrow
\B_{n,k}(\Sigma)$ (See \cite{Bir}). Similarly, the inclusion
$B_{n;k}(D)\hookrightarrow B_{n;k}(\Sigma)$ induces a monomorphism
$\B_{n;k}(D)\hookrightarrow \B_{n;k}(\Sigma)$ which will be regarded as
an inclusion.

We first determine an epimorphism $\psi_D:\B_{n;k}(D)\to H_D$ to extend
the map $\phi_D:\B_{n,k}(D)\to G_D$ for the classical
braid groups. Since we want to obtain
the homology linear representations for the classical braid groups,
we should use the $H_D=G_D$ and all of extra generators
$\bar\sigma_1,\ldots,\bar\sigma_{n-1}$ for $\B_{n;k}(D)$ should be sent
to the identity by $\psi_D$ as we have seen earlier in \S3.1.
Then $\psi_D|_{\B_{n,k}(D)}=\phi_D$. For
some extension $H$ of $G_D$, let $\psi:\B_{n;k}(\Sigma)\to H$ be an
epimorphism. In order to obtain an extension of homology linear
representations of the classical braid groups via $\psi$, the
following condition is required:
\begin{equation}
    \psi|_{\B_{n;k}(D)}=\psi_D\tag{$**$}
\end{equation}
\begin{figure}[ht]
$$
\xymatrix
{
\B_{n,k}(D)\:\ar[d]^{\phi_D}\ar@{^{(}->}[rr]&&\B_{n;k}(D)\:\ar[d]^{\psi_D}
\ar@{^{(}->}[rr]&&\B_{n;k}(\Sigma)\ar[d]^{\psi}\\
G_D\ar@{=}[rr]&&H_D\:\ar@{^{(}->}[rr]&&H
}
$$
\caption{Quotients of braid groups on a disk $D$ and surface $\Sigma$}
    \label{fig:commutativediagram}
\end{figure}

This condition is nothing but a reintepretation of Definition~\ref{defn:extension}
and is necessary to make the diagram in Figure~\ref{fig:commutativediagram}
commutative so that
$\psi|_{\B_{n,k}(D)}=\phi_D$ and we can apply the construction of
Theorem~\ref{thm:representation}.
We first
show that $(**)$ uniquely determines the quotient group $H$ for $k\ge3$.

\begin{thm}
\label{thm:quotient} Suppose $k\ge 3$. Let $\psi:\B_{n;k}(\Sigma)\to
H$ be an arbitrary epimorphism onto a
group $H$ satisfying $(**)$. Then $H$ is isomorphic to $H_\Sigma$ defined in \S3.1.
\end{thm}
\begin{proof}
Let $\psi_\Sigma:\B_{n;k}(\Sigma)\to H_\Sigma$ be the epimorphism defined in \S3.1.
We want to show $\psi:\B_{n;k}(\Sigma)\to H$ factor through $H_{\Sigma}$,
that is, there is an epimorphism $h:H_{\Sigma}\to H$ such that
$h\psi_\Sigma=\psi$.

Recall the presentation for $\B_{n;k}(\Sigma)$ in Lemma~\ref{lem:presentationofb_n;k}. The condition $(**)$ implies
$\psi(\sigma_i)=q, \psi(\zeta_j)=t$, and $\psi(\bar\sigma_m)=1$ for
all $1\le i\le k-1, 1\le j\le n$, and $1\le m\le n-1$. Since $k\ge
3$, the relation CR1 among generators in $X_2$ is not vacuous and so
the relation CR1 through CR3 for $X_2$ and the condition $(**)$ imply
$$
    [\psi(\mu_r),q]=[\psi(\lambda_r),q]=[\psi(\mu_r),t]=[\psi(\lambda_r),t]=[q,t]=1
$$
for all $1\le r\le g$. Also the relation (iii) in Lemma~\ref{lem:presentationofb_n;k} implies
$$
[\psi(\bar \mu_r), q]=[\psi(\bar \lambda_r),q]=
[\psi(\bar \mu_r), t]=[\psi(\bar \lambda_r),t]=1
$$
for all $1\le r\le g$. Thus $q$ and $t$ lie in the center of $H$.
Using this, all other relations in $H_{\Sigma}$ can be shown to hold
in $H$. Therefore $\psi$ induces an epimorphism $h:H_{\Sigma}\to H$.

To prove that $h$ is an isomorphism, it suffices to show that for any
word $W$ in generators of $H_{\Sigma}$, $h(W)=1$ implies $W=1$.

Using the relations of $H_{\Sigma}$, we may assume that $W$ is a
word of the following form:
$$
    W=q^c t^d \prod_{i=1}^g m_i^{a_i}\ell_i^{b_i}
    \bar m_i^{\bar a_i}\bar\ell_i^{\bar b_i}
$$
And this is equivalent to
$$
    m_r^{a_r}=WW_r^{-1}
    \left( \ell_r^{b_r}\bar m_r^{\bar a_r}\bar \ell_r^{\bar b_r} \right)^{-1}
$$
where $W_r=\left(q^c t^d \prod_{i\neq r} m_i^{a_i}
\ell_i^{b_i}\bar m_i^{\bar a_i}\bar \ell_i^{\bar b_i}\right)$.
The relation $\left[ m_r, \bar \ell_r \right]=q$ implies
$m_r^{a_r}\bar \ell_r=q^{a_r}\bar \ell_r m_r^{a_r}$. Then
for $1\le r\le g$,
$$
    WW_r^{-1}
\left( \ell_r^{b_r}
\bar m_r^{\bar a_r}\bar \ell_r^{\bar b_r}\right)^{-1} \bar \ell_r
=m_r^{a_r}\bar \ell_r=q^{a_r}\bar \ell_r m_r^{a_r}=q^{a_r} \bar \ell_r WW_r^{-1}
\left(
\ell_r^{\beta_r}\bar m_r^{\bar a_r}
\bar \ell_r^{\bar b_r}\right)^{-1}.
$$
Now apply $h$ and use $h(W)=1$ to have:
\begin{align*}
h(W_r^{-1})h\left( \ell_r^{b_r}\bar m_r^{\bar a_r}\bar \ell_r^{\bar b_r}\right)^{-1} h(\bar \ell_r)
&=h(q^{a_r}\bar \ell_r)h(W_r^{-1})
h\left( \ell_r^{b_r}\bar m_r^{\bar a_r}\bar \ell_r^{\bar b_r} \right)^{-1}\\
&=
h(W_r^{-1})
h(q^{a_r})
h\left( \ell_r^{b_r}\bar m_r^{\bar a_r}\bar \ell_r^{\bar b_r} \right)^{-1} h(\bar \ell_r)
\end{align*}

Thus $h(q^{a_r})=1$. Since $q$ is of infinite order by $(**)$, $a_r=0$.
Similarly, we can show that $b_r=\bar a_r=\bar b_r=0$ in
turn. Then the condition $(**)$ implies $c=d=0$. Consequently, $W=1$.
\end{proof}

For $k\le 2$, the condition $(**)$ cannot uniquely determine a quotient group of $\B_{n;k}(\Sigma)$. To take advantage of representations in analyzing the surface braid group $\B_{0,n}(\Sigma)$, one may prefer a simpler coefficient ring as long as representations carries enough information.
For the classical case, there are also several groups satisfying the condition $(*)$ unless we assume abelian. For the surface braid groups, we cannot obtain any interesting representation if an abelian coefficient ring is used as we discussed in \S2. Thus it is natural to require a kind of minimality for the choice of a quotient $H$ of $\B_{n;k}(\Sigma)$ in the sense that any further quotient of $H$ violates the condition $(**)$. Then $H_\Sigma$ given in \S3.1 is minimal for all $k$ and $n$. The proof of this minimality is similar to that of Theorem~\ref{thm:quotient}.
Consequently, Theorem~\ref{thm:quotient} forces us to use the epimorphism $\psi_\Sigma:\B_{n;k}(\Sigma)\to H_\Sigma$.

We now discuss possible actions of $\B_{0,n}(\Sigma)$ on $H_\Sigma$ induced from $\psi_\Sigma$.

\begin{thm}\label{thm:character}
Let $\psi_\Sigma:\B_{n;k}(\Sigma)\to H_\Sigma$ be the epimorphism defined in \S3.1. Let $\beta\cdot h$ denote any action on $h\in H_\Sigma$ by $\beta\in\B_{0,n}(\Sigma)$ that is induced from $\psi_\Sigma$ and satisfies the two conditions given above Theorem~\ref{thm:uniqueness}. Then
$$\beta\cdot h=h\chi(\beta)\psi_\Sigma(\beta)$$
for some function $\chi:\B_{0,n}(\Sigma)\to
C_{H_\Sigma}(G_\Sigma)$ with the property that $(\chi,\psi_\Sigma):\B_{0,n}(\Sigma)\to C_{H_\Sigma}(G_\Sigma)\rtimes H_\Sigma$ is a homomorphism
where $C_{H_\Sigma}(G_\Sigma)$ denotes the centralizer of $G_\Sigma$ in $H_\Sigma$.
\end{thm}
\begin{proof}
By the hypotheses of the action, we have
$$h'(\beta\cdot 1)=\beta\cdot (1h')=(\beta\cdot 1)\beta_\sharp(h')$$
and
$$\beta_\sharp(h')=\psi_\Sigma(\beta)^{-1} h'\psi_\Sigma(\beta)$$
for all $h'\in G_\Sigma$.
By combinig two equations, we have
$$
\psi_\Sigma(\beta)^{-1}h'\psi_\Sigma(\beta)=(\beta\cdot 1)^{-1}h'(\beta\cdot 1).
$$
and so $(\beta\cdot 1)\psi_\Sigma(\beta)^{-1}\in C_{H_\Sigma}(G_\Sigma)$.
Hence $(\beta\cdot 1)=\chi(\beta)\psi_\Sigma(\beta)$ for a function $\chi:\B_{0,n}(\Sigma)\to C_{H_\Sigma}(G_\Sigma)$.
Since $\chi(\beta_1\beta_2)\psi_\Sigma(\beta_1\beta_2)
=(\beta_1\beta_2)\cdot1
=\beta_2\cdot(\beta_1\cdot1)
=(\chi(\beta_1)\psi_\Sigma(\beta_1))\chi(\beta_2)\psi_\Sigma(\beta_2)$, we have
$$
\chi(\beta_1\beta_2)
=\chi(\beta_1)\psi_\Sigma(\beta_1)\chi(\beta_2)\psi_\Sigma(\beta_1)^{-1}.
$$
This implies that
\begin{align*}
\left(\chi(\beta_1\beta_2),\psi_\Sigma(\beta_1\beta_2)\right)&=
\left(\chi(\beta_1)\psi_\Sigma(\beta_1)\chi(\beta_2)\psi_\Sigma(\beta_1)^{-1},\psi_\Sigma(\beta_1\beta_2) \right)\\
&=\left(\chi(\beta_1),\psi_\Sigma(\beta_1)\right)
\left(\chi(\beta_2),\psi_\Sigma(\beta_2)\right).
\end{align*}
Therefore $(\chi,\psi_\Sigma):\B_{0,n}(\Sigma)\to C_{H_\Sigma}(G_\Sigma)\rtimes H_\Sigma$ is a homomorphism.
\end{proof}

The function $\chi$ in the above theorem behaves like a character of $\B_{0,n}(\Sigma)$.
In fact, if $k\ge 2$, then it can be shown that $C_{H_\Sigma}(G_\Sigma)=Z(H_\Sigma)=\langle q\rangle\oplus\langle t\rangle$. Hence $\chi$ can be any homomorphism from $\B_{0,n}(\Sigma)$ to $Z(H_\Sigma)$.
In this case, the representations $\Psi_k$ obtained from $\psi$ is given by $\Psi_k=\chi\otimes\Phi_k$ for some character $\chi$ where $\Phi_k$ is the proposed representation in Theorem~\ref{thm:representation}.

\end{document}